%% file: final_draft.tex
\documentclass[a4paper,twoside,11pt,reqno]{amsart}

\input{commands} 
\usepackage{caption}
\title{$S$-arithmetic groups acting simply transitively on products of Bruhat--Tits trees}
\author{Jonah Mendel}
\author{Jiahui Yu}
\thanks{Email addresses: Jonah Mendel, \texttt{jm316@rice.edu}; Jiahui Yu, \texttt{jiahu878@mit.edu}.}
\date{June 8, 2026}

\begin{document}

\begin{abstract}
    Lubotzky asked when a transitive $\fp$-arithmetic action on a
    Bruhat--Tits tree arising from a totally definite quaternion algebra
    can be refined to a simply transitive action. We answer this question
    affirmatively for class number one maximal orders $\fO$ in totally
    definite quaternion algebras $B$ over number fields $K$, using the
    classification of Kirschmer and Voight. For all but finitely many
    prime ideals $\fp$, we construct a congruence arithmetic lattice,
    commensurable with $\fO[1/\fp]^\times/\fo_K[1/\fp]^\times$, acting
    simply transitively on the vertices of the associated Bruhat--Tits
    tree. We further prove an $S$-arithmetic analogue for products of
    Bruhat--Tits trees, building on the work of Rungtanapirom, Stix, and
    Vdovina. These simply transitive actions produce many new examples of
    Ramanujan Cayley graphs and Ramanujan Cayley regular cubical complexes
    arising from totally definite quaternion algebras.
\end{abstract}

\maketitle

\bibliographystyle{amsplain}

\section{Introduction}

Famously, Margulis \cite{Margulis1988ExpandersConcentrators} and
Lubotzky--Phillips--Sarnak (LPS) \cite{LPS88} independently construct
explicit families of $(p+1)$-regular Ramanujan Cayley graphs for primes
$p \equiv 1 \bmod 4$, with generating sets arising from elements of norm $p$
in the Lipschitz order
\begin{equation} \label{eq:lipschitz-order}
    \fO = \bZ + \bZ\cdot i + \bZ\cdot j + \bZ\cdot ij,
    \qquad i^2 = j^2 = -1,\ ij = -ji.
\end{equation}
Central to these constructions is the fact that $\fO$ has
\emph{class number one}, which implies that, for every
odd prime $p$, the $p$-arithmetic group
\[
    \Gamma_p = \fO[1/p]^\times/\bZ[1/p]^\times
\]
acts transitively on the vertices of the Bruhat--Tits tree
$\cT(\bQ_p)$. They find that the principal congruence subgroup $\Gamma_p(2)$ acts simply
transitively on the vertices of $\cT(\bQ_p)$, and the quotients
$\Gamma_p(2q)\backslash \cT(\bQ_p)$ then give the desired Cayley graphs. The
Ramanujan property is proved using the Jacquet--Langlands correspondence and
Deligne's proof of the Ramanujan--Petersson conjecture \cite{Deligne69}.

The Margulis--LPS construction extends to other totally definite quaternion
algebras. More precisely, let $K$ be a totally real number field, let $B/K$ be a
totally definite quaternion algebra, and let $\fO \subseteq B$ be an order of class
number one. If $\fp$ is a prime ideal of $\fo_K$ such that
\[
    \fO_\fp \cong M_2(\fo_{K,\fp}),
\]
then the $\fp$-arithmetic group
\[
    \Gamma_{\fO,\fp} = \fO[1/\fp]^\times/\fo_K[1/\fp]^\times
\]
acts transitively on the vertices of the Bruhat--Tits tree $\cT(K_\fp)$.
For suitable auxiliary ideals $\fq$, the congruence quotients
\[
    \Gamma_{\fO,\fp}(\fq)\backslash \cT(K_\fp)
\]
are Ramanujan graphs by work of Jordan--Livn\'{e} \cite{JordanLivne00}
 and Blasius \cite{Blasius06}. However, they
need not carry a natural Cayley graph structure.

Toward constructing Cayley Ramanujan graphs, Lubotzky asks when one can
find a subgroup $\Gamma \leq \Gamma_{\fO,\fp}$ that acts simply
transitively on the vertices of $\cT(K_\fp)$
\cite[Remark 7.4.4]{Lubotzky94}. This question has remained largely open,
despite connections with constructions in post-quantum cryptography
\cite{CharlesGorenLauter09-Families-Of-Ramanujan-Graphs,
CharlesLauterGoren09-Cryptographic-Hash-Functions,JoYamasaki18}.

This problem becomes tractable through the
work of Kirschmer, Voight, and Lorch. In \cite{KirschmerVoight10},
Kirschmer and Voight enumerate the totally definite Eichler orders of
class number one over number fields. Kirschmer and Lorch subsequently
enumerate the totally definite Gorenstein orders of type number one or two
over number fields \cite{KirschmerLorch16}. Taken together, these results
give a complete classification of totally definite class number one
quaternion orders over number fields.

Using the Kirschmer--Voight classification, we answer Lubotzky's question
for maximal orders of class number one.

\begin{thmABC} \label{thm:cayley-ramanujan-graphs}
    Let $K$ be a totally real number field, let $B/K$ be a totally
    definite quaternion algebra, and let $\fO \subseteq B$ be a maximal
    order of class number one. Then, for all but finitely many prime
    ideals $\fp \subseteq \fo_K$, there exists a congruence arithmetic
    lattice
    \[
        \Gamma \leq \PGL_2(K_\fp)
    \]
    commensurable with $\Gamma_{\fO,\fp}$ that acts simply transitively
    on the vertices of $\cT(K_\fp)$. In particular, for each such $\fp$,
    congruence quotients of $\Gamma$ give an infinite family of
    $(N\fp+1)$-regular Ramanujan Cayley graphs whose generating sets
    come from elements $\alpha \in \fO$ with $\Nrd(\alpha)\cdot\fo_K = \fp$.
\end{thmABC}

In this paper, we also consider higher-dimensional analogues of
Ramanujan graphs, namely \emph{Ramanujan regular cubical complexes}, as
defined by Jordan--Livn\'{e} \cite{JordanLivne00} and by
Rungtanapirom, Stix, and Vdovina (RSV)
\cite{RungtanapiromStixVdovina19}. These complexes are constructed as
quotients of products of Bruhat--Tits trees by suitable
$S$-arithmetic groups, in a manner parallel to the Margulis--LPS
construction.

The construction is as follows. Let $K$, $B$, and $\fO$ be as before,
and let $S$ be a finite set of primes, with $S_0 \subseteq S$ denoting
the subset of primes at which $B$ splits. Then the $S$-arithmetic group
\[
    \Gamma_{\fO,S}
    =
    \fO[1/S]^\times/\fo_K[1/S]^\times
\]
acts transitively on the vertices of
\[
    \cT_{S_0}
    =
    \prod_{\fp \in S_0} \cT(K_\fp).
\]
By Jordan--Livn\'{e}, the quotients by principal congruence subgroups,
\[
    \cX(\fa) = \Gamma_{\fO,S}(\fa)\backslash \cT_{S_0},
\]
are Ramanujan complexes. Rungtanapirom, Stix, and Vdovina establish the
analogous Ramanujan property in the parallel global function field
setting.

To obtain a Cayley structure on $\cX(\fa)$, one seeks an intermediate
lattice
\[
    \Gamma_{\fO,S}(\fa) \leq \Lambda_S \leq \Gamma_{\fO,S}
\]
that acts simply transitively on the vertices of $\cT_{S_0}$. In this
case $\cX(\fa)$ acquires the structure of a Ramanujan \emph{Cayley}
cubical complex; see Section \ref{sec:ramanujan-cubical-complexes} for the
definition, and compare \cite[Section 3.1]{HsiehEtAl25} for a slightly
different convention.

For some time, the main well-studied examples of such Ramanujan Cayley
complexes were those arising from the Margulis--LPS generating sets
associated with the Lipschitz order. In dimension $2$ (i.e. $\#S_0 = 2$),
these complexes provide important arithmetic examples in the broader
study of lattices acting on products of two trees, as developed in work
of Burger--Mozes
\cite{BurgerMozes97,BurgerMozes00a,BurgerMozes00b} and in related work
of Burger--Mozes--Zimmer \cite{BurgerMozesZimmer09},
Kimberley--Robertson \cite{KimberleyRobertson02}, and Rattaggi
\cite{Rattaggi04}. However, one obstacle to working over $\bQ$ is that
the tree factors in $\cT_{S_0}$ have different valencies.

Strikingly, Stix and Vdovina construct a new family of explicit examples
in positive characteristic \cite{StixVdovina17}. Specifically, for each
odd prime power $q$, they construct an explicit, torsion-free
quaternionic lattice
\[
    \Lambda \leq \PGL_2(\bF_q((t))) \times \PGL_2(\bF_q((t)))
\]
that acts simply transitively on the associated product
of trees. Here, the two tree factors have the same valency. Stix and
Vdovina also apply their construction to the non-archimedean
uniformization of smooth projective surfaces over $\bF_q((t))$,
producing a non-classical fake quadric over $\bF_3((t))$.

Rungtanapirom, Stix, and Vdovina extend the Stix--Vdovina construction
to sets $S$ with arbitrarily many split places, producing torsion-free
$S$-arithmetic groups that act simply transitively on higher-dimensional
products of Bruhat--Tits trees of the same valency
\cite{RungtanapiromStixVdovina19}. 
We continue this line of work by extending Theorem \ref{thm:cayley-ramanujan-graphs} to the
higher-dimensional setting.

\begin{thmABC} \label{thm:cayley-cubical-complexes}
    Let $K$ be a totally real number field, let $B/K$ be a totally
    definite quaternion algebra, and let $\fO \subseteq B$ be a maximal
    order of class number one. Then, for every finite set of prime
    ideals $S$ disjoint from a fixed finite set of prime ideals, there
    exists a congruence arithmetic lattice
    \[
        \Gamma \leq \prod_{\fp \in S_0} \PGL_2(K_\fp)
    \]
    commensurable with $\Gamma_{\fO,S}$ that acts simply transitively
    on $\cT_{S_0}$.
\end{thmABC}

The lattices constructed by Stix--Vdovina and Rungtanapirom--Stix--Vdovina give the first 
explicit examples of irreducible arithmetic groups acting simply transitively on products of
trees of the same valency.
Our work similarly gives the first such examples for arithmetic groups defined over number fields.
Specifically, given a proper extension $K/\bQ$,
we can choose multiple distinct prime ideals $\fp$ with
isomorphic completions $K_\fp$. Applying Theorem \ref{thm:cayley-cubical-complexes} to such sets
$S$ gives simply transitive lattices acting on products of Bruhat--Tits
trees of the same valency.

\begin{thmABC} \label{thm:same-valency-products}
    For each $1 \leq n \leq 5$, there exist infinitely many
    $S$-arithmetic lattices derived from quaternion algebras over
    number fields that act simply transitively on a product of $n$
    trees of the same valency.
\end{thmABC}

The bound on $n$ comes from the fact that the
fields in the Kirschmer--Voight classification have degree at most $5$.

Our work over number fields appears to fit into the framework for
non-archimedean uniformization studied by Stix and Vdovina; see
\cite[Remark 53]{StixVdovina17}. This analogy raises natural questions
about non-archimedean uniformizations of certain complex surfaces of
small Euler characteristic. These questions, as well as torsion-free
infinite families, are pursued in joint work in preparation by the first
author with Stix and Vdovina \cite{StixMendelVdovinaInProgress}.

In another direction, Bondarenko, Grigorchuk, and Vdovina introduce
Ramanujan subshifts and use the RSV lattices to construct explicit
examples \cite{BondarenkoGrigorchukVdovina26}. Their construction makes
essential use of the equal-valency feature of the RSV lattices. The
lattices constructed in Theorem \ref{thm:same-valency-products} should provide new
number-field examples to investigate from this perspective.

Finally, in theoretical computer science, there has been growing
interest in Ramanujan cubical complexes because of their connections to
breakthroughs in error-correcting codes
\cite{DinurEtAl22a,DinurEtAl22b,DinurLinVidick24}. In one application,
Hsieh, Lubotzky, Mohanty, Reiner, and Zhang use Ramanujan Cayley
cubical complexes constructed from the Margulis--LPS generating sets to
produce the first explicit examples of constant-degree lossless vertex
expanders \cite{HsiehEtAl25}. In a related direction, Panteleev and
Kalachev propose studying sheaf codes on the RSV Ramanujan cubical complexes
as a possible approach to the
quantum locally testable code conjecture
\cite[Section VIII]{PanteleevKalachev24}.

The main technical innovation of this paper is the notion of a
\emph{complementary triple}, which generalizes features
of the constructions of Margulis, Lubotzky--Phillips--Sarnak, Rattaggi,
and Rungtanapirom--Stix--Vdovina. Its purpose is to reduce the problem
of promoting a transitive lattice action to a simply transitive one to
a finite-group computation. Theorems \ref{thm:cayley-ramanujan-graphs}, \ref{thm:cayley-cubical-complexes}, and
\ref{thm:same-valency-products} all follow from a detailed study of complementary
triples carried out in Sections \ref{sec:complementary-triples} and
\ref{sec:computing-complementary-triples}. Furthermore, drawing on the work of Rattaggi
\cite{Rattaggi04}, Rungtanapirom--Stix--Vdovina
\cite{RungtanapiromStixVdovina19}, and Chari \cite{Chari20}, we
explain how metacommutation in quaternion arithmetic yields
presentations for our lattices.

This paper is organized as follows. Section~\ref{sec:quaternion-preliminaries}
collects the necessary background on quaternion algebras and orders.
Section~\ref{sec:arithmetic-groups-products-trees} describes the natural
$S$-arithmetic action on products of Bruhat--Tits trees and explains
how orders of class number one give rise to transitive actions and
group presentations. Section~\ref{sec:ramanujan-cubical-complexes} constructs
Ramanujan regular cubical complexes from these actions.
Section~\ref{sec:complementary-triples} introduces complementary triples and shows
how they produce simply transitive lattices.
Section~\ref{sec:computing-complementary-triples} describes the finite-quotient
computations used to find complementary triples and classifies the
minimal triples associated to the maximal orders in the
Kirschmer--Voight classification. Section~\ref{sec:torsion-free-groups}
establishes torsion obstructions for lattices obtained from
complementary triples.

Our paper makes substantial use of computations carried out in Magma
\cite{BosmaCannonPlayoust97}. The Magma code used for this project is available in
the following GitHub repository:
\begin{center}
    \url{https://github.com/jonah-mendel/complementary_triples_code}.
\end{center}
Appendix~\ref{app:complementary-triple-tables}
collects the main computational output of the project.

We conclude the introduction with a concrete example illustrating the
methods of this paper. Let
\[
    B = \quat{\bQ}{-2}{-5},
\]
which has discriminant $5$ and class number one, and fix a maximal order $\fO \subseteq B$.
For any finite set $S$ of primes not containing $5$, the
$S$-arithmetic group $\Gamma_{\fO,S}$ acts transitively on the vertices
of $\cT_S$, and the vertex stabilizers are conjugate to
\[
    \fO^\times/\{\pm 1\} \cong C_3.
\]
Reduction modulo $5$ induces a map
\[
    \pi_5:\Gamma_{\fO,S} \to G[5],
    \qquad
    G[5] = (\fO/5\fO)^\times/\bF_5^\times.
\]
In Magma, we find that $\pi_5$ is injective on
$\fO^\times/\{\pm 1\}$ and that $G[5]$ contains a complementary
subgroup $H$. It follows that
\[
    \Gamma_{\fO,S}(5,H)=\pi_5^{-1}(H)
\]
acts simply transitively on the vertices of $\cT_S$. 
In Section~\ref{sec:arithmetic-groups-products-trees}, we give an explicit presentation for
$\Gamma_{\fO,S}(5,H)$ with $S = \{2,3\}$.

\begin{nota}
    We use the following notation throughout the paper.
    \begin{enumerate}
        \item For a commutative ring $A$, let $A^\times$ be the group
        of units of $A$ and let
        $\Delta(A) = A^\times/(A^\times)^2$.
        \item Let $\bF_q$ be the finite field with $q$ elements.
        \item For a CW-complex $X$, let $X^{(k)}$ be its $k$-skeleton.
        \item For a Dedekind domain $\fo$, let $\Cl(\fo)$ be its ideal
        class group.
        \item If $\fa$ is a nonzero fractional ideal of a Dedekind
        domain, let $N(\fa)$ be its ideal norm.
        \item Let $C_n$ be the cyclic group of order $n$.
        \item Let $D_n$ be the dihedral group of order $2n$.
        \item Let $S_n$ and $A_n$ be the symmetric and alternating
        groups on $n$ letters.
    \end{enumerate}
\end{nota}

\section{Quaternion preliminaries} \label{sec:quaternion-preliminaries}

In this section, we give the necessary background on quaternion algebras.
Standard references are Maclachlan--Reid \cite{MaclachlanReid03},
Vign\'{e}ras \cite{Vigneras80}, and Voight \cite{Voight21}.

\subsection{Quaternion algebras}

Let $K$ be a field. A \emph{quaternion algebra} $B$ over $K$ is a
central simple $K$-algebra of dimension $4$. 
If $\characteristic(K) \neq 2$, then $B$ admits a presentation
\[
    B = \quat{K}{a}{b}
    = K + K i + K j + K ij
\]
with $a,b \in K^\times$ and
\[
    i^2 = a,\qquad j^2 = b,\qquad ij = -ji.
\]
The algebra $B$ has a standard involution
\[
    \alpha \mapsto \overline{\alpha},
\]
called \emph{conjugation}, and the \emph{reduced trace} and
\emph{reduced norm} are
\[
    \trd(\alpha)=\alpha+\overline{\alpha},\qquad
    \Nrd(\alpha)=\alpha\overline{\alpha}.
\]
These belong to $K$ for all $\alpha \in B$. In the above basis, if
$\alpha=x+yi+zj+wij$, then
\[
    \overline{\alpha}=x-yi-zj-wij,
\]
and hence,
\[
    \trd(\alpha)=2x,\qquad
    \Nrd(\alpha)=x^2-ay^2-bz^2+abw^2.
\]
The map $\trd$ is $K$-linear, and $\Nrd$ is a quadratic form on $B$. A
field extension $L/K$ splits $B$ if $B \otimes_K L \cong M_2(L)$. For
example, $L=K(\sqrt{a})$ splits $B$, and there is an embedding
\begin{equation} \label{eq:quaternion-matrix-embedding}
    \rho : B \to M_2(L), \quad x + yi + zj + wij \mapsto
    \begin{pmatrix}
        x + y\sqrt{a} & b(z + w\sqrt{a}) \\
        z - w\sqrt{a} & x - y\sqrt{a}
    \end{pmatrix},
\end{equation}
such that $\Nrd = \det\circ\rho$ and $\trd = \tr\circ\rho$. The reduced
norm restricts to a quadratic form on
\[
    B_0 = \left\{\alpha \in B : \trd(\alpha) = 0\right\},
\]
and $B^\times$ acts on $(B_0, \Nrd|_{B_0})$ by conjugation.

\begin{prop} \label{prop:quaternion-units-special-orthogonal} \cite[Proposition 4.5.10]{Voight21}
    Let $B$ be a quaternion algebra over $K$. Then, there is a short exact sequence
    \[
        1 \to K^\times \to B^\times \to \SO(\Nrd|_{B_0})(K) \to 1.
    \]
\end{prop}

\subsection{Orders and ideals}

Let $\fo$ be a Dedekind domain whose field of fractions $K$ is either a
global or local field. A finitely generated $\fo$-module
$\fI \subseteq B$ is called an $\fo$-\emph{ideal} if
$\fI \otimes_\fo K = B$. An $\fo$-ideal $\fO \subseteq B$ is an
\emph{$\fo$-order} if it is also a unital subring. An $\fo$-order is
\emph{maximal} if it is not properly contained in another $\fo$-order,
and \emph{Eichler} if it is the intersection of two maximal
$\fo$-orders. When the context is clear, we drop the prefix $\fo$.
It follows from \cite[Section 10.3]{Voight21} that
$\Nrd(\alpha),\trd(\alpha) \in \fo$ for all $\alpha \in \fO$.

\begin{lem}\cite[Section 2.2]{MaclachlanReid03} \label{lem:maximal-orders-matrix-algebra}
    For each ideal class $[\fa] \in \Cl(\fo)$, there is a unique
    conjugacy class of maximal orders in $M_2(K)$ represented by
    \[
        \begin{pmatrix}
            \fo & \fa^{-1} \\
            \fa & \fo
        \end{pmatrix}.
    \]
    In particular, if $\fo$ is a PID, then there is one maximal order in
    $M_2(K)$ up to conjugation.
\end{lem}

Fix an $\fo$-ideal $\fI$. The \emph{left order} and \emph{right
order} of $\fI$ are respectively
\[
    \fO_L(\fI)
    =
    \left\{\alpha \in B : \alpha \fI \subseteq \fI\right\},
    \qquad
    \fO_R(\fI)
    =
    \left\{\alpha \in B : \fI\alpha \subseteq \fI\right\}.
\]
One immediately sees that $\fO_L(\fI)$ and $\fO_R(\fI)$ are
$\fo$-orders, and $\fI$ is called an
$\fO_L(\fI),\fO_R(\fI)$-ideal. If $\fJ$ is another $\fo$-ideal with
$\fO_L(\fJ) = \fO_R(\fI)$, then $\fI\cdot\fJ$ is naturally a
$\fO_L(\fI), \fO_R(\fJ)$-ideal. We give some basic
properties of ideals.
\begin{enumerate}[label=(\alph*)]
    \item $\fI$ is \emph{two-sided} if
    $\fO_L(\fI) = \fO_R(\fI)$.

    \item $\fI$ is \emph{principal} if
        \[
            \fI = \alpha \cdot \fO_R(\fI) = \fO_L(\fI) \cdot \alpha
        \]
        for some $\alpha \in B$.

    \item $\fI$ is \emph{integral} if one of the three equivalent
    properties holds \cite[Lemma 16.2.8]{Voight21}:
        \begin{enumerate}[label=(\roman*)]
            \item $\fI \subseteq \fO_L(\fI) \cap \fO_R(\fI)$;
            \item $\fI \subseteq \fO_L(\fI)$ and is a left
            $\fO_L(\fI)$-ideal in the usual sense;
            \item $\fI \subseteq \fO_R(\fI)$ and is a right
            $\fO_R(\fI)$-ideal in the usual sense.
        \end{enumerate}

    \item $\fI$ is \emph{invertible} if there exists an
    $\fO_R(\fI),\fO_L(\fI)$-ideal $\fI^{-1}$ such that
        \[
            \fI^{-1}\cdot \fI = \fO_R(\fI)
            \qquad\text{and}\qquad
            \fI\cdot \fI^{-1} = \fO_L(\fI).
        \]
\end{enumerate}
Fortunately, when $\fO_L(\fI)$ or $\fO_R(\fI)$ is maximal, the following
lemma tells us that invertibility is automatic.

\begin{lem} \cite[Lemma 16.6.15]{Voight21} \label{lem:maximal-left-right-order-invertibility}
    Let $\fI$ be an $\fo$-ideal. If $\fO_L(\fI)$ or $\fO_R(\fI)$ is
    maximal, then both are maximal and $\fI$ is invertible.
\end{lem}

The \emph{reduced norm} of an ideal is the fractional $\fo$-ideal
generated by the reduced norms of its elements,
\[
    \Nrd(\fI)
    =
    \sum_{\alpha \in \fI} \Nrd(\alpha)\cdot \fo.
\]
If $\fI$ is integral, then $\Nrd(\fI) \subseteq \fo$. If
$\fI,\fJ$ are invertible ideals with
$\fO_L(\fJ) = \fO_R(\fI)$, then
$\Nrd(\fI\fJ) = \Nrd(\fI)\Nrd(\fJ)$. For proofs, see
\cite[Lemmas 16.3.2, 16.3.7, 16.6.15]{Voight21}.

\subsection{Quaternions over local fields} \label{sec:quaternions-local-fields}

We now consider quaternion algebras over local fields.
For archimedean fields, the only non-split case is
\[
    B \cong \quat{\bR}{-1}{-1},
\]
and from Proposition \ref{prop:quaternion-units-special-orthogonal}, we see that
$B^\times/\bR^\times \cong \SO(3)$ is compact.

Now let $\fo$ be a complete discrete valuation ring with valuation
$v:K^\times \to \bZ$, field of fractions $K$, maximal ideal
$\fp=\varpi \fo$, and residue field $\fo/\fp$.

\subsubsection{The split case and the Bruhat--Tits tree}
\label{sec:split-local-case}

Suppose $B \cong M_2(K)$ is split, and set $\fO_0=M_2(\fo)$.
Then, $M_2(K)$ has infinitely many maximal orders, all conjugate to
$\fO_0$ by Lemma~\ref{lem:maximal-orders-matrix-algebra}.

The integral right $\fO_0$-ideals of reduced norm $\fp^e$ are of the
form
\begin{equation} \label{eq:split-local-right-ideal-normal-form}
    \begin{pmatrix}
        \varpi^u & c \\
        0 & \varpi^v
    \end{pmatrix} \fO_0,
    \quad u,v \in \bZ_{\geq 0}, \quad u+v=e,
\end{equation}
where $c$ ranges over representatives of $\fo/\varpi^u \fo$. There are
exactly
\[
    N(\fp)^e + N(\fp)^{e-1} + \cdots + N(\fp) + 1
\]
such ideals.

If $\fI$ is an integral right $\fO_0$-ideal, then
$\fO_L(\fI)$ is a maximal order by Lemma~\ref{lem:maximal-left-right-order-invertibility}.
Conversely, every maximal order in $M_2(K)$ arises in this way.
Multiplying $\fI$ by a scalar does not change $\fO_L(\fI)$, so for each
maximal order $\fO$ we may choose such an ideal $\fI$ with
$\fI \not\subseteq \varpi\fO_0$. The ideal $\Nrd(\fI)=\fp^e$
is independent of the choice of $\fI$, 
and we say that
$\fO_0 \cap \fO$ is an Eichler order of \emph{level} $\fp^e$.

The \emph{Bruhat--Tits tree} $\cT(K)$ is the graph whose vertices are
maximal orders in $M_2(K)$. Two vertices $\fO$ and $\fL$ are joined by
an edge if and only if $\fO \cap \fL$ is an Eichler order of level
$\fp$. Equivalently, $\fO$ and $\fL$ are adjacent if and only if there
exists an integral $\fO,\fL$-ideal $\fP$ with $\Nrd(\fP)=\fp$.

More generally, the distance from $\fO_0$ to $\fL$ is the unique
integer $e$ such that $\fO_0 \cap \fL$ is an Eichler order of level
$\fp^e$.

\begin{thm}[{\cite[Proposition 23.5.8]{Voight21}}]
    The graph $\cT(K)$ is a tree. Moreover, $\PGL_2(K)$ acts
    transitively on $\cT(K)$ by conjugation, and the stabilizer of a
    vertex $\fO$ is $\fO^\times/\fo^\times$.
\end{thm}

A standard reference on the Bruhat--Tits tree is Serre \cite{Serre80}.

\subsubsection{The non-split case} \label{sec:nonsplit-local-case}

Suppose $B$ is a division quaternion algebra over $K$. Then the
valuation on $K$ determines a unique integer-valued valuation on $B$,
which we normalize by
\[
    v_B(\alpha)=v(\Nrd(\alpha)) \in \bZ.
\]
The unique maximal order of $B$ is
\[
    \fO = \left\{\alpha \in B : v_B(\alpha) \geq 0\right\}.
\]
Moreover, $B$ has a unique two-sided ideal of reduced norm $\fp$, given
by
\[
    \fP = \left\{\alpha \in B : v_B(\alpha) > 0\right\}.
\]
In particular, the only Eichler order in $B$ is the maximal order. One
can explicitly write
\[
    B = \left(\frac{u,\varpi}{K}\right),
\]
where $K(u)/K$ is the unique unramified quadratic extension, and
\[
    \fP = j\fO = \fO j, \quad \fP^{2} = \fp\fO.
\]

We now inspect the structure of the unit groups $B^\times$ and
$\fO^\times$. First, it follows from Proposition \ref{prop:quaternion-units-special-orthogonal} that
$B^\times/K^\times$ is compact. For positive integers $k$, the
\emph{higher $k$-unit groups} of $\fO$ are $1 + \fP^k$. Set
\[
    U^k_\fp = \fO^\times/(\fo^\times\cdot (1 + \fP^{2k})),
    \quad
    \wtU^k_\fp = B^\times/(K^\times \cdot (1 + \fP^{2k})).
\]
For each $k$, we have a split short exact sequence
\begin{equation} \label{eq:ramified-projective-unit-extension}
    1 \to U^k_\fp \to \wtU^k_\fp \xrightarrow{\bar v_B} C_2 \to 1,
\end{equation}
where $\bar v_B(\alpha) \equiv v_B(\alpha) \bmod 2$. A section of
$\bar v_B$ is given by sending the nontrivial element of $C_2$ to $j$.

Observe that the groups $U_\fp^k$ and $\wtU_\fp^k$ are completely
determined by the choice of $\fp$ and $k$, since $B/K$ is unique. Both
groups have a construction similar to the finite matrix groups
$\PGL_2(\fo/\fp^k)$. However, unlike $\PGL_2(\fo/\fp^k)$, both groups
are always solvable.

\begin{prop} \label{prop:ramified-projective-units-solvable}
    The groups $U_\fp^k$ and $\wtU_\fp^k$ are solvable and
    \[
        \# U_\fp^{k} = N\fp^{3k - 1}(N\fp+1).
    \]
\end{prop}

\begin{proof} 
    Since $[\wtU_\fp^k:U_\fp^k] = 2$, we need only show that $U_\fp^k$ is solvable. 
    For positive integers $k$, we have short exact sequences 
    \[ 
        1 \to \frac{1+\fP^{k-1}}{1+\fP^k} \to \frac{\fO^\times}{1+\fP^k} \to \frac{\fO^\times}{1+\fP^{k-1}} \to 1. 
    \] 
    Then 
    \[ 
        \frac{\fO^\times}{1+\fP} \cong \bF_{q^2}^\times \quad\text{and}\quad \frac{1+\fP^{k-1}}{1+\fP^k} \cong \bF_{q^2} 
    \] 
    where $q = N\fp$. 
    The solvability and order of $U_\fp^k$ now follow from induction on $k$. 
\end{proof}

\begin{cor}
    For $k \geq l$, the kernel of $\wtU_\fp^k \to \wtU_\fp^l$ is a
    $p$-group, where $p$ is the characteristic of $\fo/\fp$.
\end{cor}

We have $\fo/\fp^k \subseteq \fO/\fP^{2k}$, so
$\fO/\fP^{2k}$ is a $\fo/\fp^k$-algebra. Let
\[
    G_\fO[\fp^k] = (\fO/\fP^{2k})^\times/(\fo/\fp^k)^\times.
\]
When the context is clear, we drop the $\fO$. The following lemma
identifies $U_\fp^k$ with $G_\fO[\fp^k]$.

\begin{lem} \label{lem:ramified-reduction-isomorphism}
    The map
    \[
        \eta: (\fO/\fP^{2k})^\times
        \to \fO^\times/(1 + \fP^{2k}),
        \quad
        \alpha + \fP^{2k} \mapsto \alpha (1 + \fP^{2k})
    \]
    is an isomorphism which induces an isomorphism
    \[
        \bar\eta: G[\fp^k] \longrightarrow U_\fp^k.
    \]
\end{lem}

\begin{proof}
    Since $\fP$ is the Jacobson radical of $\fO$, an element of $\fO$
    is a unit if and only if its image in $\fO/\fP^{2k}$ is a unit.
    Thus, $\eta$ is a well-defined surjection, and its kernel is
    $1+\fP^{2k}$. Therefore, $\eta$ is an isomorphism. Since
    $\fo \cap \fP^{2k}=\fp^k$, the map $\eta$ sends
    $(\fo/\fp^k)^\times$ to $\fo^\times\cdot(1+\fP^{2k})$, so it
    induces the claimed isomorphism $\bar\eta$.
\end{proof}

The reduced norm descends to the groups $U_\fp^k$ and $G_\fO[\fp^k]$.
Indeed, for $\alpha \in \fO$, observe that
\[
    \Nrd(1 + \varpi^k\alpha)
    =
    (1 + \varpi^k\alpha)(1 + \varpi^k\overline\alpha)
    =
    1 + \varpi^k\trd(\alpha) + \varpi^{2k}\Nrd(\alpha).
\]
Therefore, $\Nrd(1 + \fP^{2k}) \subseteq 1 + \fp^k$, and $\Nrd$
descends to a map
\[
    \overline{\Nrd}:  U_\fp^k \to \Delta(\fo/\fp^k).
\]
On the other hand, let $\omega_1,\ldots,\omega_4$ be an $\fo$-basis for
$\fO$ (which exists because $\fo$ is a PID). Then the reduced norm on
$\fO$ can be written as the quadratic form
\[
    Q_\fO(X_1,X_2,X_3,X_4)
    =
    \sum_i \Nrd(\omega_i)X_i^2
    +
    \sum_{j < k} \trd(\omega_j\overline{\omega_k})X_jX_k.
\]
Since $\omega_i \in \fO$, the coefficients of $Q_\fO$ are in $\fo$, and
by reducing coefficients modulo $\fp^k$, we get a quadratic form
$\overline{Q}_\fO$ on $\fO/\fP^{2k}$. From the above discussion, we get
the following commutative triangle.
\[\begin{tikzcd}
	{G[\fp^k]} && {U_\fp^k} \\
	& {\Delta(\fo/\fp^k)} & {}
	\arrow["{\bar\eta}", from=1-1, to=1-3]
	\arrow["{\overline{Q}_{\fO}}"', from=1-1, to=2-2]
	\arrow["{\overline{\Nrd}}", from=1-3, to=2-2]
\end{tikzcd}\]
We further investigate the forms $\overline{Q}_\fO$ and
$\overline{\Nrd}$ in Section \ref{sec:computing-complementary-triples}.

\subsection{Quaternions over number fields} \label{sec:quaternions-number-fields}

Let $K$ be a number field, let $V(K)$ be the set of places of $K$, and set
\begin{align*}
    V_\infty(K) &= \left\{v \in V(K) : \text{$v$ is archimedean}\right\}, \\
    V_f(K) &= \left\{v \in V(K) : \text{$v$ is finite}\right\}.
\end{align*}
If every archimedean place of $K$ is real, then $K$ is \emph{totally
real}. The \emph{discriminant} of $K$ is denoted $d_K$.
For a finite set $S \subseteq V_f(K)$, the ring of \emph{$S$-integers}
of $K$ is
\[
    \fo_K[1/S]
    =
    \left\{x\in K : v_\fp(x)\geq 0 \text{ for all } \fp \notin S\right\}.
\]
We write $\fo_K=\fo_K[1/\emptyset]$, which is the usual ring of integers.

For $v \in V(K)$, let $B_v = B \otimes_K K_v$. A quaternion algebra
$B/K$ is \emph{split} at $v$ if $B_v$ splits, and \emph{ramified}
otherwise. Let $\Ram(B)$ be the set of ramified places of $B$, and set
\begin{align*}
    \Ram_\infty(B) &= \Ram(B) \cap V_\infty(K), \\
    \Ram_f(B) &= \Ram(B) \cap V_f(K).
\end{align*}
It follows from class field theory that $B \mapsto \Ram(B)$ induces a
bijection
\begin{equation} \label{eq:ramification-classification}
    \left\{
        \begin{array}{c}
        \text{isomorphism classes of} \\
        \text{quaternion algebras } B/K
        \end{array}
    \right\}
    \longleftrightarrow
    \left\{
        \begin{array}{c}
        \text{finite subsets of non-complex places of } K \\
        \text{of even cardinality}
        \end{array}
    \right\}.
\end{equation}
See \cite[Theorem 14.6.1]{Voight21} for details. The
\emph{discriminant} of $B$ is the integral $K$-ideal
\[
    \fD(B) = \prod_{\fp \in \Ram_f(B)} \fp.
\]
For a finite set $S\subseteq V_f(K)$, we say that $B$ is
\emph{$S$-definite} if
\[
    S\cup V_\infty(K) \subseteq \Ram(B)
\]
and \emph{$S$-indefinite} otherwise. If
$\Ram_\infty(B) = V_\infty(K)$, then we say that $B$ is
\emph{totally definite}. If $B$ is totally definite, then $K$ is
totally real, since a complex place cannot be in $\Ram_\infty(B)$.
It follows from \eqref{eq:ramification-classification} that a totally definite $B$ is
completely determined by $\fD(B)$ up to isomorphism. We also write
$D(B) = N(\fD(B))$.

Let $\fo = \fo_K[1/S]$ for some finite $S\subseteq V_f(K)$, and let
$\fI$ be an $\fo$-ideal. We write
$\fI_\fp = \fI \otimes_\fo \fo_{K,\fp}$ for $\fp \notin S$. The
\emph{local-global principle} tells us that
\[
    \fI = \bigcap_{\fp \not\in S} \fI_\fp.
\]
In particular, an $\fo$-order $\fO$ is maximal if and only if
$\fO_\fp$ is maximal for all $\fp \notin S$. Also, an arbitrary
$\fo$-order is maximal at all but finitely many places. For more
details on the local-global principle, see
\cite[Theorem 9.4.9]{Voight21}.

Let $\fE = \fO \cap \fO'$ be an Eichler $\fo$-order. Then
$\fE_\fp$ is Eichler for each $\fp \notin S$, and $\fE_\fp$ is maximal
for $\fp \in \Ram_f(B) \setminus S$. The \emph{level} of $\fE$ is
\[
    \fN(\fE) = \prod_{\fp \notin S} \fN(\fE_\fp).
\]
Because $\fE_\fp$ is maximal for all but finitely many $\fp$, the level
$\fN(\fE)$ is an integral $\fo$-ideal. Furthermore, $\fN(\fE)$ and
$\fD(B)$ are coprime, since $\fE_\fp$ is maximal for all
$\fp \mid \fD(B)\cdot \fo_K[1/S]$.

\begin{lem} \label{lem:local-left-order-invertibility}
    Let $\fI$ be an integral right $\fO$-ideal. If
    $\fO_R(\fI)_\fp$ is maximal for all $\fp \mid \Nrd(\fI)$, then
    $\fI$ is invertible.
\end{lem}

\begin{proof}
    It suffices to check invertibility locally. If
    $\fp \nmid \Nrd(\fI)$, then $\fI_\fp=\fO_R(\fI)_\fp$. If
    $\fp \mid \Nrd(\fI)$, then $\fO_R(\fI)_\fp$ is maximal by
    assumption, so $\fI_\fp$ is invertible by Lemma~\ref{lem:maximal-left-right-order-invertibility}.
    Hence, $\fI$ is locally invertible at every finite place, and
    therefore invertible.
\end{proof}

\subsection{The ideal class number and type number}

We briefly review the ideal class number and type number of an order. For details,
see \cite[Chapters 17, 18, 27]{Voight21}.

Let $K$ be a number field and let $\fo = \fo_K[1/S]$ be the
$S$-integers of $K$ for some finite set $S \subseteq V_f(K)$. Let
$B/K$ be a quaternion algebra, and let $\fO \subseteq B$ be an
$\fo$-order. The \emph{genus} of $\fO$ is the set $\cG(\fO)$ of all
$\fo$-orders $\fL \subseteq B$ such that, for every nonzero prime
ideal $\fp$ of $\fo$, the localizations $\fL_\fp$ and $\fO_\fp$ are
isomorphic as $\fo_\fp$-algebras. For example, if $\fO$ is an Eichler
order of level $\fN$, then $\cG(\fO)$ consists of all Eichler orders in
$B$ of level $\fN$. The group $B^\times$ acts on $\cG(\fO)$ by
conjugation, and the \emph{type number} of $\fO$ is
\[
    t(\fO) = \#(B^\times \backslash \cG(\fO)),
\]
which is finite \cite[Corollary 27.6.25]{Voight21}.

Let $\cR(\fO)$ denote the set of invertible right $\fO$-ideals. Then
$B^\times$ acts on $\cR(\fO)$ by left multiplication, and the
\emph{class number} of $\fO$ is
\[
    h(\fO) = \#(B^\times \backslash \cR(\fO)),
\]
which is also finite \cite[Main Theorem 27.6.14]{Voight21}. There is a
natural $B^\times$-equivariant surjective map
\begin{equation} \label{eq:right-ideal-to-left-order-map}
    \cR(\fO) \to \cG(\fO), \quad \fI \mapsto \fO_L(\fI).
\end{equation}
In particular, $t(\fO) \leq h(\fO)$, and $h(\fO) = h(\fL)$ if
$\fL \in \cG(\fO)$.

\section{Arithmetic groups acting on products of trees} \label{sec:arithmetic-groups-products-trees}

From now on, we fix a totally real number field $K$,
a totally definite quaternion algebra $B/K$,
and a finite set $S \subseteq V_f(K)$ with
\[
    S_0 = S \setminus \Ram_f(B), \qquad S_1 = S \cap \Ram_f(B).
\]
The goal of this section is to describe $S$-arithmetic groups acting
simply transitively on
\[
    \cT_{S_0}=\prod_{\fp \in S_0}\cT(K_\fp).
\]

\subsection{Unit groups and normalizers}

Every $\fo_K[1/S]$-order is of the form
\[
    \fO[1/S]
    = \fO \otimes_{\fo_K} \fo_K[1/S]
\]
for an $\fo_K$-order $\fO$. We say that an $\fo_K$-order $\fO$ is
\emph{$S$-maximal} if $\fO_\fp$ is maximal for all $\fp \in S$. If
$\fO$ is $S$-maximal, then it follows from
Lemma~\ref{lem:local-left-order-invertibility} that every integral right
$\fO$-ideal $\fI$ with $\Nrd(\fI)$ supported on $S$ is invertible.

Fix an $S$-maximal order $\fO$, and let
\[
    N_{\fO,S} = N_{B^\times}(\fO[1/S])/K^\times, \qquad
    \Gamma_{\fO,S} = \fO[1/S]^\times/\fo_K[1/S]^\times, \qquad
    \Gamma_{\fO,S}^1 = \fO[1/S]^1/\{\pm 1\}.
\]
When the context is clear, we drop the subscript $\fO$. We have a map,
not necessarily surjective,
\[
    \Nrd:\Gamma_S \longrightarrow \Delta^+(\fo_K[1/S]),
\]
with kernel $\Gamma_S^1$, where
\[
    \Delta^+(\fo_K[1/S])
    =
    \left\{[x] \in \Delta(\fo_K[1/S]) :
    \text{$v(x) > 0$ for all $v \in V_\infty(K)$}\right\}.
\]
Thus, $\Gamma_S^1$ is a finite index normal subgroup in $\Gamma_S$. We
refer the reader to \cite[Section 32.2]{Voight21} for details.

The group $\Gamma_S$ is also normal and finite index in $N_S$, but the
quotient is slightly more involved for general orders $\fO$. However,
we have a clean description for Eichler orders of class number one.

\begin{prop}\label{prop:normalizer-unit-exact-sequence}
    Let $\fO[1/S]$ be an Eichler $\fo_K[1/S]$-order of class number one
    and level $\fN$. Then there is a short exact sequence
    \begin{equation}\label{eq:normalizer-unit-exact-sequence}
        1 \longrightarrow \Gamma_S \longrightarrow N_S \longrightarrow
        \prod_{\fp \mid \fD(B)\,\fN\, \fo_{K,S}} C_2 \longrightarrow 1.
    \end{equation}
    In particular, if $\fO[1/S]$ is maximal and
    $S_1=\Ram_f(B)$, then $\Gamma_S \cong N_S$.
\end{prop}

\begin{proof}
    By \cite[Proposition 18.5.3]{Voight21}, the fact that
    $\fO[1/S]$ has class number one gives a short exact sequence
    \[
        1 \to \Gamma_S \to N_S
        \to \Pic_{\fo_K[1/S]}(\fO[1/S]) \to 1,
    \]
    where $\Pic_{\fo_K[1/S]}(\fO[1/S])$ is the Picard group. By
    \cite[23.4.19]{Voight21}, this Picard group is the desired product
    of cyclic groups of order $2$.
\end{proof}

Assume the hypotheses of Proposition~\ref{prop:normalizer-unit-exact-sequence}. For
later use, suppose that $T \subseteq S$ and
$S\setminus T \subseteq \Ram_f(B)$. Write
$S\setminus T=\{\fp_1,\dots,\fp_r\}$, and let
$\fP_i=\alpha_i\fO$ be the corresponding two-sided $\fO$-ideal of
reduced norm $\fp_i$. Applying \eqref{eq:normalizer-unit-exact-sequence} with $T$ in
place of $S$, the image of each $\alpha_i$ is nontrivial in exactly one
$C_2$-factor, and these factors are distinct for different $i$.
Consequently,
\[
    \Gamma_S
    =
    \left\langle
        \Gamma_T \cup \left\{\alpha_i : \fp_i \in S\setminus T\right\}
    \right\rangle.
\]

By Borel--Harish-Chandra \cite{BorelHarishChandra62}, the diagonal
embedding realizes $N_S$ as a cocompact, discrete subgroup in
\[
    \prod_{v \in S\cup V_\infty(K)} B_v^\times/K_v^\times.
\]
Since $B$ is totally definite, the factors over $V_\infty(K)$ are
compact. When $v \in S_1$, the factors $B_v^\times/K_v^\times$ are also
compact. Projecting away from the compact places in
$S_1\cup V_\infty(K)$, we realize $N_S$ as a discrete, cocompact subgroup in
\[
    \prod_{\fp \in S_0} \PGL_2(K_\fp)
\]
acting properly and cocompactly on $\cT_{S_0}$.

We want a global description of the vertices of $\cT_{S_0}$. For this
purpose, let
\[
    \cG_S(\fO) = \left\{\fL \in \cG(\fO) : \fL[1/S] = \fO[1/S]\right\}.
\]
By the local-global principle, there is an $N_S$-equivariant bijection
\[
    \cG_S(\fO) \to \cT_{S_0}^{(0)}, \quad
    \fL \mapsto (\fL_\fp)_{\fp \in S_0}.
\]
For the rest of the paper, we identify $\cT_{S_0}^{(0)}$ with
$\cG_S(\fO)$.

\subsection{Simply transitive actions and integral ideals}
We first establish some notation. 
For an ideal $\fI$ of $B$ and $\alpha \in B^\times$, let
\[
    \fI^\alpha = \alpha\fI\alpha^{-1}.
\]

Let $\fO$ be an $S$-maximal order as before.
For $\fp \in S_0$, we call $\fL \in \cG_S(\fO)$ a
\emph{$\fp$-neighbor} of $\fO$ if there exists an
$\fL,\fO$-integral ideal $\fP$ with $\Nrd(\fP) = \fp$. Thus, $\fO$ and
$\fL$ are adjacent in $\cT_{S_0}$ if and only if $\fL$ is a
$\fp$-neighbor of $\fO$ for some $\fp \in S_0$.
We now turn our attention to the action of $\Gamma_S$ on
$\cG_S(\fO)$. Our exposition here is heavily influenced by
Stix--Vdovina \cite{StixVdovina17},
Rungtanapirom--Stix--Vdovina \cite{RungtanapiromStixVdovina19},
and Kirschmer--Voight \cite{KirschmerVoight10}.

First, we determine when the action is transitive.

\begin{prop} \label{prop:class-number-one-transitivity}
    Let $B$ be a totally definite quaternion algebra, let
    $S\subseteq V_f(K)$ be finite, and let $\fO \subseteq B$ be an
    $S$-maximal order. For the action of $\Gamma_S$ on
    $\cG_S(\fO)$, the stabilizer of $\fO$ is $\Gamma_{S_1}$.
    Moreover, $\Gamma_S$ acts transitively on $\cG_S(\fO)$ if and only
    if every integral right $\fO[1/S_1]$-ideal $\fP$ with
    $\Nrd(\fP) = \fp \cdot \fo_K[1/S_1]$ for some $\fp \in S_0$ is
    principal.
\end{prop}

\begin{proof}
    We first compute the stabilizer of $\fO$. The stabilizer of $\fO$
    in $N_S$ is $N_{S_1}$, so the stabilizer of $\fO$ in $\Gamma_S$ is
    \[
        N_{S_1} \cap \Gamma_S.
    \]
    We claim that this is precisely $\Gamma_{S_1}$. Let
    $\alpha \in N_{S_1} \cap \Gamma_S$. Then, for every
    $\fp \in S_0$, we have
    \[
        \fO_\fp^\alpha=\fO_\fp.
    \]
    Thus, $\alpha$ normalizes $\fO_\fp$. After multiplying $\alpha$ by
    a scalar in $\fo_K[1/S]^\times$, we may assume that
    $\alpha \in \fO_\fp^\times$ for all $\fp \in S_0$.
    Since $\alpha \in \fO[1/S]$, we also have $\alpha \in \fO_\fp$ for
    all $\fp \notin S$, which implies that
    \[
        \alpha \in \bigcap_{\fp \notin S_1} \fO_\fp = \fO[1/S_1].
    \]
    Therefore, the class $\alpha\cdot \fo_K[1/S]^\times$ lies in
    $\Gamma_{S_1}$. The reverse inclusion is immediate.

    We now prove the transitivity criterion. Suppose first that every
    integral right $\fO[1/S_1]$-ideal $\fI$ with
    $\Nrd(\fI)=\fp \cdot \fo_K[1/S_1]$ for some $\fp \in S_0$ is
    principal. We begin with the distance one case. Let $\fL$ be a
    $\fp$-neighbor of $\fO$ for some $\fp \in S_0$. Then
    $(\fO\cap\fL)[1/S_1]$ is an Eichler $\fo_K[1/S_1]$-order. By the
    construction in Section \ref{sec:quaternions-number-fields}, there is a unique
    integral $\fo_K[1/S_1]$-ideal $\fI$ such that
    \begin{equation} \label{eq:neighbor-ideal-local-conditions}
        \fO_R(\fI) = \fO[1/S_1], \quad
        \fO_L(\fI) = \fL[1/S_1], \quad
        \Nrd(\fI) = \fp\cdot\fo_K[1/S_1].
    \end{equation}
    By assumption, $\fI$ is principal, say
    \[
        \fI = \alpha \cdot \fO[1/S_1].
    \]
    Since the left order of $\fI$ is $\fL[1/S_1]$, we have
    \[
        \fO[1/S_1]^\alpha = \fL[1/S_1].
    \]
    Hence, $\fO^\alpha=\fL$, so every $\fp$-neighbor of $\fO$ lies in
    the $\Gamma_S$-orbit of $\fO$.

    Now let $\fL \in \cG_S(\fO)$ be arbitrary. Since $\cT_{S_0}$ is
    connected, choose a path
    \[
        \fO=\fO_0,\fO_1,\ldots,\fO_r=\fL
    \]
    in $\cT_{S_0}$. We prove by induction on $i$ that $\fO_i$ lies in
    the $\Gamma_S$-orbit of $\fO$. The case $i=0$ is clear. Suppose
    that $\fO_i=\fO^\gamma$ for some $\gamma \in \Gamma_S$. Since
    $\fO_{i+1}$ is adjacent to $\fO_i$, the order
    $\fO_{i+1}^{\gamma^{-1}}$ is adjacent to $\fO$. By the distance one
    case, there exists $\alpha \in \Gamma_S$ such that
    \[
        \fO^\alpha = \fO_{i+1}^{\gamma^{-1}}.
    \]
    We obtain $\fO_{i+1}\in \Gamma_S\cdot \fO$. Thus, by induction,
    $\fL$ lies in the $\Gamma_S$-orbit of $\fO$. Therefore, $\Gamma_S$
    acts transitively on $\cG_S(\fO)$.

    Conversely, suppose that $\Gamma_S$ acts transitively on
    $\cG_S(\fO)$. Fix $\fp \in S_0$. Let $\fL$ be a $\fp$-neighbor of
    $\fO$. By transitivity, there exists
    $\alpha \in \fO[1/S]^\times$ such that
    \[
        \fO^\alpha=\fL.
    \]
    After multiplying $\alpha$ by a scalar in $\fo_K[1/S]^\times$, we
    may assume that $\alpha \in \fO[1/S_1]$ and that
    \[
        \fI=\alpha\cdot\fO[1/S_1]
    \]
    satisfies the conditions in \eqref{eq:neighbor-ideal-local-conditions}. In
    particular, $\fI$ is a principal integral right
    $\fO[1/S_1]$-ideal of reduced norm $\fp$.

    We claim that distinct $\fp$-neighbors give distinct principal
    ideals in this way. Let $\fO_1$ and $\fO_2$ be two
    $\fp$-neighbors of $\fO$, with corresponding elements
    $\alpha_1$ and $\alpha_2$. If
    \[
        \alpha_1\cdot\fO[1/S_1]=\alpha_2\cdot\fO[1/S_1],
    \]
    then $\alpha_1=\alpha_2\beta$ for some
    $\beta \in \fO[1/S_1]^\times$. By the stabilizer computation above,
    $\beta$ fixes $\fO$. Hence,
    \[
        \fO_1=\fO^{\alpha_1}
        =\fO^{\alpha_2\beta}
        =\fO^{\alpha_2}
        =\fO_2,
    \]
    a contradiction. Thus, distinct $\fp$-neighbors give distinct
    principal integral right $\fO[1/S_1]$-ideals of reduced norm $\fp$.

    Since there are $N\fp+1$ such neighbors, we have constructed
    $N\fp+1$ distinct principal integral right $\fO[1/S_1]$-ideals of
    reduced norm $\fp$. By the local-global principle and the ideal
    count in Section \ref{sec:split-local-case}, these are all such integral
    ideals. Therefore, every integral right $\fO[1/S_1]$-ideal $\fP$
    with $\Nrd(\fP)=\fp\cdot\fo_K[1/S_1]$ for some $\fp \in S_0$ is
    principal.
\end{proof}

\begin{cor}
    If $\fO[1/S_1]$ is class number one, then $\Gamma_{\fO,S}$ acts
    transitively on $\cG_S(\fO)$.
\end{cor}

Let $\Gamma \leq \Gamma_{\fO,S}$ act simply transitively on
$\cG_S(\fO)$. Once such a subgroup exists, the ideal arithmetic of
$\fO$ determines its generators and relations. For $\fp \in S_0$, let
\[
    A_\fp=\left\{\text{integral right $\fO[1/S_1]$-ideals of reduced norm $\fp$}\right\},
    \quad
    A = \bigcup_{\fp \in S_0} A_\fp.
\]
Then $\#A_\fp=1+N\fp$ by Section \ref{sec:split-local-case}. There is a
bijection between $A_\fp$ and the $\fp$-neighbors of $\fO$ given by
\[
    \fP \longmapsto \fO_L(\fP).
\]
Let $\fP \in A_\fp$. Then $\alpha \in \fO[1/S]^\times$ conjugates
$\fO$ to $\fO_L(\fP)$ if and only if
$\alpha\cdot \fO[1/S_1] = \fP$. Since $\Gamma$ acts transitively on
$\cG_S(\fO)$, it follows from Proposition \ref{prop:class-number-one-transitivity} that,
for each $\fp \in S_0$ and $\fP \in A_\fp$, there exists
$\alpha_{\fP} \in \fO[1/S_1]$ such that
\[
    \fP = \alpha_{\fP}\cdot\fO[1/S_1]
    \quad\text{and}\quad
    \alpha_{\fP}\cdot\fo_K[1/S]^\times \in \Gamma.
\]
Conversely, if $\fP = \alpha \cdot \fO[1/S_1]$ and
$\alpha\cdot\fo_K[1/S]^\times \in \Gamma$, then
$\alpha\alpha_\fP^{-1}$ stabilizes $\fO$. Since $\Gamma$ acts simply
transitively, $\alpha\alpha_\fP^{-1} \in \fo_K[1/S]^\times$. In
particular, there is a well-defined map
\[
    A_\fp \to \Gamma, \quad
    \fP \mapsto \alpha_{\fP} \cdot \fo_K[1/S]^\times.
\]
We now establish some basic facts about $\Gamma$.

\begin{lem}
    $\Gamma$ is generated by $\{\alpha_{\fP} : \fP \in A\}$.
\end{lem}

\begin{proof}
    Let $\Lambda \leq \Gamma$ be the subgroup generated by
    $\{\alpha_{\fP} : \fP \in A\}$. By construction, $\Lambda$ sends
    $\fO$ to each of its neighbors in $\cT_{S_0}$. Thus, $\Lambda$ acts
    transitively on $\cG_S(\fO)$. If $\fO^\lambda = \fO^\gamma$ for
    $\lambda \in \Lambda$ and $\gamma \in \Gamma$, then
    $\gamma^{-1}\lambda$ stabilizes $\fO$, and is thus trivial.
    Therefore, $\Lambda = \Gamma$.
\end{proof}

\begin{lem}
    Let $\fP \in A_\fp$. Then there exists a unique ideal
    $\iota(\fP) \in A_\fp$ such that
    $\alpha_{\fP}\alpha_{\iota(\fP)} \in \fo_K[1/S]^\times$.
    Furthermore, $\Gamma$ inverts the edge corresponding to $\fP$ if
    and only if $\iota(\fP) = \fP$. In this case, the stabilizer in
    $\Gamma$ of the corresponding unoriented edge is
    $\lag \alpha_{\fP} \rag \cong C_2$.
\end{lem}

\begin{proof}
    Existence follows because $\alpha_{\fP}^{-1}$ carries $\fO$ to a
    $\fp$-neighbor of $\fO$, so by simple transitivity it is represented
    by a unique element $\alpha_{\fP'}$ with $\fP'\in A_\fp$.
    We see that $\Gamma$ inverts the edge associated to $\fP$ if and
    only if $\alpha_{\fP}^2 \in \fo_K[1/S]^\times$, which gives the
    second statement.
\end{proof}

Now let $\tau$ be a two-dimensional face adjacent to $\fO$. This face
lies in the product of two distinct tree factors, say those
corresponding to $\fp,\fq \in S_0$. Choose $\fP \in A_\fp$ and
$\fQ \in A_\fq$ for the two edges of $\tau$ adjacent to $\fO$. There
are then unique ideals $\fP' \in A_\fp$ and $\fQ' \in A_\fq$ which
describe the opposite $\fp$- and $\fq$-edges of the face. More
precisely, the opposite edges are represented by
$(\fP')^{\alpha_{\fQ}}$ and $(\fQ')^{\alpha_{\fP}}$. The situation is
summarized in Figure~\ref{fig:tree-product-square}.

\begin{figure}[htbp]
    \centering
    \begin{tikzpicture}[
        x=3.0cm,
        y=2.0cm,
        >=stealth,
        edge/.style={->, shorten >=5pt, shorten <=5pt, line width=0.65pt, draw=blue!65!black},
        vertex/.style={circle, draw=black, fill=orange!75, inner sep=0pt, minimum size=5pt},
        edge label/.style={font=\small, fill=white, inner sep=1.5pt},
        vertex label/.style={font=\small, inner sep=2pt}
        ]

        \coordinate (Ocoord) at (0,0);
        \coordinate (Pcoord) at (1,0);
        \coordinate (Qcoord) at (0,1);
        \coordinate (Tcoord) at (1,1);

        \draw[edge] (Ocoord) -- node[edge label, below=4pt] {$\fP$} (Pcoord);
        \draw[edge] (Ocoord) -- node[edge label, left=4pt] {$\fQ$} (Qcoord);
        \draw[edge] (Qcoord) -- node[edge label, above=4pt] {$(\fP')^{\alpha_{\fQ}}$} (Tcoord);
        \draw[edge] (Pcoord) -- node[edge label, right=4pt] {$(\fQ')^{\alpha_{\fP}}$} (Tcoord);

        \node[vertex] at (Ocoord) {};
        \node[vertex] at (Pcoord) {};
        \node[vertex] at (Qcoord) {};
        \node[vertex] at (Tcoord) {};

        \node[vertex label, below left=2pt] at (Ocoord) {$\fO$};
        \node[vertex label, below right=2pt] at (Pcoord) {$\fO^{\alpha_{\fP}}$};
        \node[vertex label, above left=2pt] at (Qcoord) {$\fO^{\alpha_{\fQ}}$};
        \node[vertex label, above right=2pt] at (Tcoord) {$\fO_\tau$};

        \node[font=\small, fill=white, inner sep=1pt] at (0.5,0.5) {$\tau$};

    \end{tikzpicture}
    \caption{A face of $\cT_{S_0}$ adjacent to $\fO$, with one edge in
    the $\fp$-direction and one edge in the $\fq$-direction.}
    \label{fig:tree-product-square}
\end{figure}

Let $\fO_\tau$ be the vertex of $\tau$ opposite $\fO$. The two paths
around the face give
\[
    \fO_\tau=\fO^{\alpha_{\fQ}\alpha_{\fP'}}
    =\fO^{\alpha_{\fP}\alpha_{\fQ'}}.
\]
Since the two products $\alpha_{\fP}\alpha_{\fQ'}$ and
$\alpha_{\fQ}\alpha_{\fP'}$ both send $\fO$ to $\fO_\tau$ and
$\Gamma$ acts simply transitively on vertices, they differ by a scalar
in $\fo_K[1/S]^\times$. Hence,
\begin{equation} \label{eq:square-face-scalar-relation}
    \alpha_{\fP}\alpha_{\fQ'}\alpha_{\fP'}^{-1}\alpha_{\fQ}^{-1}
    \in \fo_K[1/S]^\times.
\end{equation}
Furthermore, given $\alpha_{\fP}$ and $\alpha_{\fQ}$, it follows from
Chari's work on metacommutation \cite[Section 3]{Chari20} that
\[
    \fP' = \alpha_{\iota(\fQ)}\fP + \fp \fO, \quad
    \fQ' = \alpha_{\iota(\fP)}\fQ + \fq \fO.
\]
Combining the previous lemmas with \cite[Theorem 1]{Brown84}, we can
compute a presentation of $\Gamma$.

\begin{thm} \label{thm:simply-transitive-presentation}
    The map $\fP \mapsto \alpha_{\fP}$ induces an isomorphism
    \[
        \left\langle
            A
            \ \middle|\
            \begin{array}{l}
                \fP\cdot \iota(\fP), \\
                \fP\cdot \fQ' \cdot \iota(\fP') \cdot \iota(\fQ)
                \text{ satisfying \eqref{eq:square-face-scalar-relation}}
            \end{array}
        \right\rangle
        \cong \Gamma.
    \]
\end{thm}

\begin{rmk}
    In Theorem \ref{thm:simply-transitive-presentation}, we treat the ideals as formal
    letters. The relations do not describe a product of compatible
    ideals. However, we can write the product of compatible ideals
    \[
        (\fQ')^{\alpha_\fP} \cdot \fP =
        (\fP')^{\alpha_\fQ} \cdot \fQ.
    \]
\end{rmk}

\begin{rmk}
    Groups that act simply transitively on a product of two trees are
    sometimes called Burger--Mozes--Wise (BMW) groups. For more on BMW
    groups, see the survey by Caprace \cite{Caprace19}.
\end{rmk}

We now give some examples. In describing our examples, we borrow
terminology introduced in Section \ref{sec:complementary-triples}.

\begin{ex} \label{ex:lipschitz-presentation}
    We return to the example from the introduction. Let
    $B = \quat{\bQ}{-2}{-5}$. Then $B$ is a class number one
    quaternion algebra with discriminant $5$. Fixing a maximal order
    $\fO$, we find in Magma that $\Gamma_\emptyset \cong C_3$, and
    there is a complementary triple $(\emptyset,5,H)$ of $\fO$. Setting
    $S = \{2,3\}$, we compute
    \[
        \Gamma_{\fO,S}(5,H) \cong
        \left\langle
        a_1,a_2,a_3,b_1,b_2,b_3,b_4
        \ \middle|\
        \begin{array}{@{}l@{\qquad}l@{\qquad}l@{}}
        a_i^2,\ 1\leq i\leq 3
        &
        a_2b_3a_1b_3
        &
        a_1b_4a_1b_2
        \\
        b_j^2,\ 1\leq j\leq 4
        &
        a_2b_1a_2b_4
        &
        a_3b_2a_2b_2
        \\
        &
        a_3b_1a_1b_1
        &
        a_3b_4a_3b_3
        \end{array}
        \right\rangle
    \]
    where $A_2 = \{a_i : 1 \leq i \leq 3\}$ and
    $A_3 = \{b_j : 1 \leq j \leq 4\}$. We can write the generators
    explicitly as
    \[
        \begin{array}{@{}ll@{\qquad\qquad}ll@{}}
        a_1 = & [i],
        &
        b_1 = & [4+4i-2j-k],
        \\[1mm]
        a_2 = & [4i-2j+k],
        &
        b_2 = & [-2+2i+k],
        \\[1mm]
        a_3 = & [4-2j-k],
        &
        b_3 = & [2j-k],
        \\[1mm]
        &
        &
        b_4 = & [2-k].
        \end{array}
    \]
    Note that we scaled some elements by a power of two to kill
    denominators.
\end{ex}

\begin{ex} \label{ex:qsqrt2-torsion-free}
    Let $K = \bQ(\sqrt 2)$ and let $B/K$ be the quaternion algebra with
    \[
        \fD(B) = \fp_2\fp_9,
    \]
    where $\fp_2$ and $\fp_9$ are the primes of norm $2$ and $9$,
    respectively. Then, for the maximal order $\fO$, we check in Magma
    that $\Gamma_\emptyset \cong C_3$, and there is a complementary
    triple $(\emptyset,\fp_2,H)$ of $\fO$. Let
    $S = \{\fp_7,\fq_7\}$ be the primes of norm $7$. Then we compute
    \[
        \Gamma_{\fO,S}(\fp_2,H) \cong
        \left\langle
        a_1,\ldots,a_8,b_1,\ldots,b_8
        \ \middle|\
        \begin{array}{llllll}
        a_1a_3
        &
        a_2a_8
        &
        a_4b_4a_3b_2
        &
        a_2b_8a_8b_2
        &
        a_7b_5a_8b_1
        &
        a_7b_6a_2b_7
        \\
        a_4a_5
        &
        a_6a_7
        &
        a_4b_1a_2b_4
        &
        a_7b_4a_6b_3
        &
        a_4b_7a_6b_8
        &
        a_4b_3a_2b_1
        \\
        b_1b_7
        &
        b_2b_8
        &
        a_1b_5a_7b_2
        &
        a_1b_2a_6b_5
        &
        a_4b_2a_3b_3
        &
        a_4b_5a_5b_6
        \\
        b_3b_4
        &
        b_5b_6
        &
        a_1b_4a_8b_6
        &
        a_1b_1a_3b_7
        &
        a_1b_6a_8b_3
        &
        a_4b_8a_7b_7
        \end{array}
        \right\rangle
    \]
    where $A_{\fp_7} = \{a_i : 1 \leq i \leq 8\}$ and
    $A_{\fq_7} = \{b_j : 1 \leq j \leq 8\}$.

    Recall that any torsion element of $\Gamma_{\fO,S}(\fp_2,H)$ must,
    after conjugacy, either invert an edge or invert a face in
    $\cT_S$ adjacent to our base vertex $\fO$ (see
    Lemma \ref{lem:simply-transitive-torsion-order-two}). The edge-inverting candidates are the
    elements $a_i$ and $b_j$, and the face-inverting candidates are the
    elements $a_i b_j$. A direct check of the relations shows that none
    of these elements has finite order. Hence,
    $\Gamma_{\fO,S}(\fp_2,H)$ is torsion-free.
\end{ex}

For function-field analogues of these constructions, Stix and Vdovina
give rank-two examples \cite{StixVdovina17}, while Rungtanapirom,
Stix, and Vdovina compute explicit presentations for the
higher-dimensional groups $\Gamma_{\fO,S}$
\cite{RungtanapiromStixVdovina19}. In their setting, $B$ is the
quaternion algebra over $\bF_q(t)$ ramified at $0$ and $1$, $\fO$ is
the maximal order over the ring of $\{0,1\}$-integers in $B$, and $S$
is a set of $\bF_q$-rational places
containing $\{0,1\}$. 

It would be interesting to compute a similar
family of presentations for the groups $\Gamma_{\fO,S}$ arising from a
fixed maximal, or more generally Eichler, order in a totally definite
quaternion algebra over a number field. We do not pursue this question
here.

\section{Ramanujan regular cubical complexes} \label{sec:ramanujan-cubical-complexes}

In this section, we briefly review the construction of Ramanujan Cayley
graphs and, more generally, of Ramanujan Cayley regular cubical
complexes.

Let $X$ be a connected $k$-regular graph with adjacency matrix $A_X$.
Recall that $X$ is \emph{Ramanujan} if every eigenvalue $\lambda$ of
$A_X$ satisfies either
\[
    \lambda = \pm k
    \quad\text{or}\quad
    |\lambda| \leq 2\sqrt{k-1}.
\]
Moreover, $X$ is bipartite if and only if $-k$ is an eigenvalue of
$A_X$.

The universal cover of a connected $k$-regular graph is the
$k$-regular tree. More generally, a cubical complex whose universal
cover is a product of regular trees is called a \emph{regular cubical
complex}. Jordan--Livn\'e \cite{JordanLivne00} and
Rungtanapirom--Stix--Vdovina \cite{RungtanapiromStixVdovina19} extend
the notion of a Ramanujan graph to this setting. More precisely,
Jordan--Livn\'e define Ramanujan local systems on regular cubical
complexes, although we will not need that level of generality here. The
approach we take to Ramanujan cubical complexes is closer to that of
Rungtanapirom--Stix--Vdovina in
\cite[Section 6]{RungtanapiromStixVdovina19}.

Suppose that $\cX$ is a finite regular cubical complex with universal
cover
\[
    \cT_{r_1} \times \cdots \times \cT_{r_m},
\]
where $\cT_{r_i}$ is the $r_i$-regular tree. Since the factors are
labeled, each edge of $\cX$ has a well-defined \emph{direction} $i$. We
say that two vertices $x,y \in \cX^{(0)}$ are \emph{$i$-adjacent},
denoted $x \sim_i y$, if they are connected by an edge of direction
$i$. For each $i$, define
\[
    \Adj_i:L^2(\cX^{(0)}) \to L^2(\cX^{(0)}),
    \qquad
    \Adj_i(f)(x) = \sum_{y \sim_i x} f(y).
\]
The sum is taken with multiplicity. That is, if there are several edges
of direction $i$ from $x$ to $y$, then $f(y)$ appears with the
corresponding multiplicity. We say that $\cX$ is \emph{Ramanujan} if,
for each $i$, every eigenvalue $\lambda$ of $\Adj_i$ satisfies either
\[
    \lambda = \pm r_i
    \quad\text{or}\quad
    |\lambda| \leq 2\sqrt{r_i-1}.
\]

We now describe the construction of Ramanujan regular cubical complexes.
This material is well-known to experts. Standard accounts in the graph
case include Lubotzky \cite{Lubotzky94} and Charles--Goren--Lauter
\cite{CharlesGorenLauter09-Families-Of-Ramanujan-Graphs}; the
corresponding construction for regular cubical complexes can be found
in Jordan--Livn\'e \cite{JordanLivne00} and in
Rungtanapirom--Stix--Vdovina
\cite[Section 6]{RungtanapiromStixVdovina19}.

Let $K$ be a number field, let $B/K$ be an $S$-indefinite quaternion algebra, and let
$\fO \subseteq B$ be an $S$-maximal $\fo_K$-order for
$S \subseteq V_f(K)$ finite. Let
$\Lambda_S \leq \Gamma_S$ be a congruence subgroup of level $\fC$
acting simply transitively on $\cG_S(\fO)$, and let $A$ be as in the
previous section. Fix an ideal $\fa \subseteq \fo_K$ coprime to $\fC$
and to the primes in $S$, and suppose that
\[
    \fO_\fq \cong M_2(\fo_{K,\fq})
    \quad\text{for each } \fq \mid \fa.
\]
Reduction modulo $\fa$ induces a map
\[
    \pi_\fa:\Gamma_S \to \PGL_2(\fo_K/\fa).
\]
We set
\[
    \Lambda_S(\fa)=\Lambda_S \cap \ker(\pi_\fa).
\]
Since $B$ is $S$-indefinite, strong approximation gives
\[
    \PSL_2(\fo_K/\fa)
    \subseteq \pi_\fa(\Lambda_S)
    \subseteq \PGL_2(\fo_K/\fa).
\]
Let
\[
    \cX = \Lambda_S\backslash \cT_{S_0},
    \qquad
    \cX(\fa) = \Lambda_S(\fa)\backslash \cT_{S_0}.
\]
Then $\cX(\fa)$ covers the one-vertex cube (orbi)-complex $\cX$. In
particular, the $1$-skeleton can be written as
\[
    \cX(\fa)^{(1)}
    =
    \Cay\left(\Lambda_S(\fa)\backslash \Lambda_S, A(\fa) \right),
    \qquad A(\fa) = \pi_\fa(A).
\]
When $\#S_0=1$, the product $\cT_{S_0}$ is a single tree and
$\cX(\fa)$ is a Ramanujan Cayley graph.

When $\#S_0 > 1$, the partition of $A(\fa)$ into directions
$\fp \in S_0$, given by
\[
    A_\fp(\fa) = \pi_\fa(A_\fp),
\]
records the cubical structure. In particular, $\cX(\fa)$ is a
\emph{Cayley} regular cubical complex in the following sense: for
$\fp_1,\fp_2 \in S_0$, the corresponding sets of generators satisfy
\begin{equation} \label{eq:metacommutation}
    A_{\fp_1}(\fa) \cdot A_{\fp_2}(\fa)
    =
    A_{\fp_2}(\fa) \cdot A_{\fp_1}(\fa).
\end{equation}
This is the notion of a Cayley regular cubical complex described in
\cite[Remark 6.16]{RungtanapiromStixVdovina19}. See also
\cite[Section 3]{HsiehEtAl25}, which uses a slightly different
formulation; our construction can be modified to produce Cayley regular
cubical complexes in their sense as well. Property
\eqref{eq:metacommutation} is sometimes called
\emph{metacommutation}, and it follows from
Theorem \ref{thm:simply-transitive-presentation}; see also Chari
\cite[Section 3]{Chari20}.

Jordan and Livn\'e show that these congruence quotients are Ramanujan,
assuming the relevant cases of the Ramanujan--Petersson conjecture for
Hilbert modular forms \cite[Theorem 3.1]{JordanLivne00}. The necessary
cases are known by work of Deligne \cite{Deligne69}, Carayol
\cite{Carayol86}, and Blasius \cite{Blasius06}.

There is an analogous construction when $K$ is a global function field.
In this setting, Rungtanapirom, Stix, and Vdovina show that the
complexes $\cX(\fa)$ are Ramanujan
\cite[Section 6]{RungtanapiromStixVdovina19}, using Drinfeld's proof
of the Ramanujan--Petersson conjecture for $\GL_2(K)$
\cite{Drinfeld88}.

To determine the quotient group and the resulting Cayley structure, we
now specialize to prime congruence levels. Let $\fq \notin S$ be a
prime of norm $q$ such that $\fO_\fq \cong M_2(\fo_{K,\fq})$.
Choose representatives
\[
    \varepsilon_1,\dots,\varepsilon_r \in \fo_K[1/S]^\times
\]
for an $\bF_2$-basis of $\Delta^+(\fo_K[1/S])$, and set
\[
    L = K(\sqrt{\varepsilon_1},\dots,\sqrt{\varepsilon_r}).
\]
It follows from Dirichlet's $S$-unit theorem
\cite[Proposition VI.1.1]{Neukirch99} that $r \geq 1$, so $L/K$ is a
nontrivial elementary abelian $2$-extension. For a prime $\fq$
unramified in $L$, let
\[
    \Frob_\fq \in \Gal(L/K)
\]
denote the Frobenius element associated to $\fq$. Reduction modulo
$\fq$ induces a map
\[
    \delta_\fq:\Delta^+(\fo_K[1/S]) \to \Delta(\bF_q)\cong C_2.
\]
Then $\delta_\fq(\varepsilon_i)$ is trivial in $\Delta(\bF_q)$ if and
only if
\[
    \Frob_\fq(\sqrt{\varepsilon_i}) = \sqrt{\varepsilon_i}.
\]
Since $B$ is $S$-indefinite, strong approximation implies that the
reduced norm map
\[
    \Nrd:\Lambda_S \to \Delta^+(\fo_K[1/S])
\]
is surjective. 

Combining the above discussion with the \v{C}ebotarev
density theorem \cite[VII.13.4]{Neukirch99} gives the following
criteria for the structure of $\pi_\fq(\Lambda_S)$.

\begin{thm} \label{thm:ramanujan-quotient-cayley-structure}
    With notation as above, the quotient $\pi_\fq(\Lambda_S)$ is
    isomorphic to $\PSL_2(\bF_q)$ if and only if $\Frob_\fq$ is trivial,
    which occurs with Dirichlet density $1/2^r$.
\end{thm}

\begin{rmk}
    In \cite{JoYamasaki18}, Jo and Yamasaki give a systematic
    construction of Cayley graphs of $\PSL_2(\bF_q)$ for suitable
    primes $q$, with generating sets arising from definite quaternion
    algebras over $\bQ$ of class number one. They call these graphs
    \emph{LPS-type graphs}. Their work is motivated in part by Cayley
    hash functions in the sense of Charles, Lauter, and Goren
    \cite{CharlesLauterGoren09-Cryptographic-Hash-Functions}. Aside
    from the known cases when $\fD(B) = 2,13$, resolved by the work of
    Margulis \cite{Margulis1988ExpandersConcentrators},
    Lubotzky--Phillips--Sarnak \cite{LPS88}, and Chiu \cite{Chiu92},
    they leave open the question of when the resulting graphs are
    Ramanujan. Combining our classification of complementary triples in
    Appendix~\ref{app:complementary-triple-tables} with
    Theorem~\ref{thm:ramanujan-quotient-cayley-structure}, we obtain infinite
    families of Ramanujan LPS-type graphs for each definite class
    number one quaternion algebra over $\bQ$.
\end{rmk}

\section{Complementary triples} \label{sec:complementary-triples}

As before, let $K$ be a number field and
$B$ a totally definite quaternion algebra.
In this section, we explain how to find congruence subgroups
$\Gamma \leq \Gamma_{\fO,S}$ which act simply transitively on products
of Bruhat--Tits trees.

\subsection{Abstract complementary triples}

Let $G$ be a finite group and let $T \subseteq \Ram_f(B)$. Suppose
that $\fO$ is a $T$-maximal order of $B$, and
there exists an injective homomorphism
\[
    \pi:\Gamma_{\fO,T} \hookrightarrow G
\]
and a finite set $\kappa(\pi) \subseteq V_f(K)$ disjoint from $T$ such that
\begin{enumerate}
    \item $\kappa(\pi)$ contains all primes $\fp$ for which $\fO_\fp$ is not maximal;
    \item for every finite set
    $S \subseteq V_f(K)\setminus \kappa(\pi)$ satisfying $S_1 = T$,
    the map $\pi$ extends to a homomorphism
    \[
        \pi_S:\Gamma_{\fO,S} \to G;
    \]
    \item there exists a subgroup $H \leq G$ \emph{complementary} to
    $\pi(\Gamma_{\fO,T})$; that is,
    \[
        H \cap \pi(\Gamma_{\fO,T}) = 1, \qquad
        H \cdot \pi(\Gamma_{\fO,T}) = G.
    \]
\end{enumerate}
We call $(T,\pi,H)$ a \emph{complementary triple} of $\fO$. We call the
triple \emph{normal} if $H$ is normal in $G$. We call two triples
$(T,\pi,H_1)$ and $(T,\pi,H_2)$ \emph{equivalent} if $H_1$ is conjugate
to $H_2$ in $G$.

Finite sets $S \subseteq V_f(K)\setminus \kappa(\pi)$ with $S_1=T$ are
called \emph{$\pi$-eligible}. We show that $\pi_S^{-1}(H)$ acts simply
transitively on $\cG_S(\fO)$ whenever $\Gamma_S$ acts transitively on
$\cG_S(\fO)$.
\begin{lem}
    Let $G$ be a group and let $H,U \leq G$ be complementary
    subgroups. Then $H \cap gUg^{-1} = 1$ for all $g \in G$.
\end{lem}

\begin{proof}
    Since $H$ and $U$ are complementary, $H$ acts transitively by left
    multiplication on $G/U$. The stabilizer of $U$ is $H\cap U=1$, so
    the action is simply transitive. The stabilizer of $gU$ is
    $H\cap gUg^{-1}$, and therefore $H\cap gUg^{-1}=1$ for all
    $g\in G$.
\end{proof}

\begin{prop} \label{prop:complementary-triple-simple-transitive}
    Let $(T,\pi,H)$ be a complementary triple of $\fO$. If $S$ is
    $\pi$-eligible and $\Gamma_S$ acts transitively on $\cG_S(\fO)$,
    then
    \[
        \Gamma_{\fO,S}(\pi,H) = \pi_S^{-1}(H)
    \]
    acts simply transitively on $\cG_S(\fO)$. Furthermore, if $H$ is
    normal in $G$, then
    \[
        \Gamma_{\fO,S}
        =
        \Gamma_{\fO,S}(\pi,H) \rtimes \Gamma_{\fO,T}.
    \]
\end{prop}

\begin{proof}
    By Proposition \ref{prop:class-number-one-transitivity}, the stabilizer of any order
    in $\cG_S(\fO)$ in $\Gamma_S$ is conjugate to $\Gamma_T$. Since $H$
    is complementary to $\pi(\Gamma_T)$, the subgroup
    $\Gamma_{\fO,S}(\pi,H)$ is complementary to $\Gamma_T$ in
    $\Gamma_S$. By the previous lemma, $\Gamma_{\fO,S}(\pi,H)$ has
    trivial intersection with all conjugates of $\Gamma_T$, so it acts
    freely on $\cG_S(\fO)$.

    Let $\fL \in \cG_S(\fO)$. Since $\Gamma_S$ acts transitively, there
    exists $\alpha \in \Gamma_S$ with $\fO^\alpha = \fL$. By
    our assumptions, there exists $\gamma \in \Gamma_T$ such that
    $\alpha\gamma \in \Gamma_{\fO,S}(\pi,H)$. Since $\gamma$
    stabilizes $\fO$, we have
    \[
        \fO^{\alpha\gamma} = \fL.
    \]
    Thus, $\Gamma_{\fO,S}(\pi,H)$ acts transitively. Since it also acts
    freely, it acts simply transitively.

    If $H$ is normal in $G$, then $\Gamma_{\fO,S}(\pi,H)$ is normal in
    $\Gamma_{\fO,S}$. Since it is complementary to
    $\Gamma_{\fO,T}$, the stated semidirect product decomposition
    follows.
\end{proof}

\begin{rmk}
    A similar construction appears in \cite[Proposition 5.19]{EvraEtAl26},
    where the authors construct level-$H$ congruence $p$-arithmetic
    lattices in unitary groups that act simply transitively on the
    \emph{hyperspecial} vertices of certain Bruhat--Tits buildings.
\end{rmk}

\subsection{Complementary ideals}

The first examples of complementary triples come from congruence quotients away
from $T$. Let
\[
    \fC = \fp_1^{k_1}\cdots \fp_r^{k_r} \subseteq \fo_K
\]
be an ideal coprime to every prime in $T$, and set
\[
    G_\fO[\fC] = (\fO/\fC \fO)^\times/(\fo_K/\fC)^\times.
\]
When the situation is clear, we drop the $\fO$. By the Sunzi remainder theorem, we have
\begin{equation} \label{eq:sunzi-factorization}
    G[\fC] \cong \prod_{i=1}^r G[\fp_i^{k_i}].
\end{equation}
If $\fC = \fp^k$ and $\fO_\fp$ is maximal, then the natural embedding
$\fO \to \fO_\fp$ induces an isomorphism
\begin{equation} \label{eq:prime-power-local-identification}
    G_\fO[\fp^k] \cong
    \begin{cases}
        \PGL_2(\fo_K/\fp^k) & \text{$B$ splits at $\fp$}, \\
        G_{\fO_\fp}[\fp^k] \cong U_\fp^k & \text{$B$ ramifies at $\fp$}.
    \end{cases}
\end{equation}
For general $\fC$, reduction modulo $\fC$ induces a map
\[
    \pi_{\fC} = \pi_{\fO,T,\fC}: \Gamma_{T} \longrightarrow G[\fC].
\]
Suppose that $\pi_\fC$ is injective and there exists a complement
$H \leq G[\fC]$ to $\pi_\fC(\Gamma_T)$. Then $(T,\pi_\fC,H)$ is a
complementary triple of $\fO[1/T]$ with
\[
    \kappa(\pi_\fC)
    =
    (\Ram_f(B) \setminus T)
    \cup
    \left\{\fp \in V_f(K): \fp \mid \fC\right\}.
\]
Indeed, if $S$ is disjoint from $\kappa(\pi_\fC)$ and satisfies
$S_1=T$, then reduction modulo $\fC$ extends $\pi_\fC$ to a homomorphism
\[
    \pi_{\fC,S}:\Gamma_{\fO,S}\to G[\fC].
\]
We call $\fC$ an \emph{(uninverted) complementary ideal of $\fO$ at $T$},
denote by $(T,\fC,H)$ the complementary triple, and call
$\pi_\fC$-eligible sets \emph{$\fC$-eligible}.

\subsection{Inverted complementary ideals}

The second source of complementary triples comes from ramified places in
$T$. Suppose that $\fC = \fp^k$ with $\fp$ ramified in $B$ and
$\fO_\fp$ maximal. Let $U_\fp^k$ and $\wtU_\fp^k$ be the quotient groups
defined in Section \ref{sec:nonsplit-local-case}. The embedding
$B \to B_\fp$ induces a homomorphism
\[
    \rho_\fC = \rho_{\fO,\fC}: \Gamma_{\Ram_f(B)}
    \longrightarrow \wtU_\fp^k.
\]
Furthermore,
\[
    \rho_{\fC}(\Gamma_{\Ram_f(B) \setminus \{\fp\}}) \subseteq U_\fp^k.
\]
The following commutative diagram shows how these maps fit together.
\[\begin{tikzcd}
	& {\Gamma_T} & \\
	{\Gamma_{\Ram_f(B) \setminus \{\fp\}}} & {\Gamma_{\Ram_f(B)}} \\
	{U_\fp^k} & {\wtU_\fp^k} & {C_2} \\
	{G[\fp^k]}
	\arrow[dashed, tail, from=1-2, to=2-1]
	\arrow[tail, from=1-2, to=2-2]
	\arrow[tail, from=2-1, to=2-2]
	\arrow["{\rho_\fC}", from=2-1, to=3-1]
	\arrow["{\pi_{\fC}}"', curve={height=18pt}, from=2-1, to=4-1]
	\arrow["{\rho_\fC}", from=2-2, to=3-2]
	\arrow[tail, from=3-1, to=3-2]
	\arrow["{\bar\eta}", from=3-1, to=4-1]
			\arrow["{\bar v_{B_\fp}}", from=3-2, to=3-3]
	\arrow[curve={height=-12pt}, from=3-3, to=3-2]
\end{tikzcd}\]
The dotted arrow exists when $\fp \notin T$. When $\fp \in T$, suppose
that $\rho_\fC$ is injective on $\Gamma_T$ and that there exists a
complement $H \leq \wtU_\fp^k$ to $\rho_\fC(\Gamma_T)$. Then
$(T,\rho_\fC,H)$ is a complementary triple of $\fO$, with
\[
    \kappa(\rho_\fC)=\Ram_f(B)\setminus T.
\]
We call $\fC$ an \emph{inverted complementary ideal of $\fO$ at $T$},
and write $(T,\fC,H)$ to denote the triple.

The next proposition shows that, under a few standard assumptions, a
complement inside $U_\fp^k$ automatically lifts to a complement inside
$\wtU_\fp^k$.

\begin{prop} \label{prop:ramified-complement-lift}
    Let $S \subseteq \Ram_f(B)$ with $\fp \in S$, and set
    $T = S \setminus \{\fp\}$. Let $\fO$ be an $S$-maximal order such
    that $\fO[1/T]$ is Eichler and class number one. Suppose that
    $\rho_\fC$ is injective on $\Gamma_T$. If
    $H \leq U_\fp^k$ is a complement of $\rho_\fC(\Gamma_T)$, then
    $\rho_\fC$ is injective on $\Gamma_S$, and $H$ is a complement of
    $\rho_\fC(\Gamma_S)$ in $\wtU_\fp^k$.

    In particular, if $\fC$ is an uninverted complementary ideal of $\fO$
    at $T$, then it is an inverted complementary ideal of
    $\fO$ at $S$.
\end{prop}

\begin{proof}
    Let $\fP$ be the unique two-sided $\fO$-ideal with
    $\Nrd(\fP)=\fp$. Since $\fO[1/T]$ has class number one, we can
    write
    \[
        \fP = \alpha \cdot \fO[1/T].
    \]
    The image $\rho_\fC(\alpha)$ has nontrivial image under
    $\bar v_{B_\fp}:\wtU_\fp^k \to C_2$, while
    $\rho_\fC(\Gamma_T) \subseteq U_\fp^k=\ker(\bar v_{B_\fp})$. By the
    discussion after Proposition \ref{prop:normalizer-unit-exact-sequence}, the group
    $\Gamma_S$ is generated by $\Gamma_T$ and $\alpha$, and
    $[\Gamma_S:\Gamma_T]=2$. It follows that $\rho_\fC$ is injective on
    $\Gamma_S$.

    Since $H$ is a complement of $\rho_\fC(\Gamma_T)$ in $U_\fp^k$, we
    have
    \[
        H \cap \rho_\fC(\Gamma_T)=1,
        \qquad
        \#H\cdot \#\rho_\fC(\Gamma_T)=\#U_\fp^k.
    \]
    Since $H\subseteq U_\fp^k$ and
    $\rho_\fC(\Gamma_S)\setminus \rho_\fC(\Gamma_T)$ lies in the
    nontrivial coset of $U_\fp^k$ in $\wtU_\fp^k$, we also have
    \[
        H \cap \rho_\fC(\Gamma_S)=1.
    \]
    Finally,
    \[
        \#H\cdot \#\rho_\fC(\Gamma_S)
        =
        2\#H\cdot \#\rho_\fC(\Gamma_T)
        =
        2\#U_\fp^k
        =
        \#\wtU_\fp^k.
    \]
    Thus, $H$ is a complement of $\rho_\fC(\Gamma_S)$ in
    $\wtU_\fp^k$.
\end{proof}

As mentioned in the introduction, Lubotzky asks in
\cite[Remark 7.4.4]{Lubotzky94} whether, whenever $\Gamma_S$ acts
transitively on $\cG_S(\fO)$, there exists a subgroup of
$\Gamma_S$ acting simply transitively on $\cG_S(\fO)$.
Lubotzky suggests that the ramified primes of $B$ are a natural
place to look for such simply transitive congruence groups. As can be
seen in Appendix \ref{app:complementary-triple-tables}, ramified
primes do indeed yield many complementary triples. This phenomenon is
not accidental.

Recall from Proposition \ref{prop:ramified-projective-units-solvable} that the group $\wtU_\fp^k$ is solvable. Hall's subgroup theorem
implies that, if $\#\rho_\fC(\Gamma_S)$ and
$[\wtU_\fp^k:\rho_\fC(\Gamma_S)]$ are coprime, then there exists a
complementary subgroup $H \subseteq \wtU_\fp^k$. Combining the above
discussion with the order of $\wtU_\fp^k$ given in
Proposition \ref{prop:ramified-projective-units-solvable}, we get the following proposition.

\begin{prop}
    Let $S$, $\fp$, and $\fO$ be as in
    Proposition \ref{prop:ramified-complement-lift}. Suppose that $\rho_\fC$ is
    injective on $\Gamma_S$. If
    \[
        \#\Gamma_S \quad\text{and}\quad
        2N\fp^{3k-1}(N\fp + 1)/\#\Gamma_S
    \]
    are coprime integers, then there exists an inverted complementary
    triple $(S,\fp^k,H)$ of $\fO$.
\end{prop}

\section{Computing complementary triples} \label{sec:computing-complementary-triples}

In this section, we describe how to compute complementary triples.
Throughout, $K$ is a number field, $B/K$ is a totally definite
quaternion algebra, and $\fO\subseteq B$ is an $\fo_K$-order. We fix a
finite set
\[
    T\subseteq \Ram_f(B)
\]
such that $\fO$ is $T$-maximal. Suppose there exists an $\fo_K$-basis
\begin{equation} \label{eq:order-basis}
    \fO=\fo_K\omega_1\oplus\cdots\oplus \fo_K\omega_4
\end{equation}
and write
\[
    \omega_i\omega_j=\sum_{\ell=1}^4 c_{ij\ell}\omega_\ell,
    \qquad c_{ij\ell}\in \fo_K.
\]
Note that such a basis always exists when $\fo_K$ is a principal ideal
domain; in particular, this holds when $\fO$ has class number one.

For an ideal $\fC\subseteq \fo_K$ coprime to every finite prime in $T$,
we explain how to compute the reduction modulo $\fC$ map
\[
    \pi_{\fC}:\Gamma_T\longrightarrow G[\fC].
\]
Then, for $\fC = \fp^k$ with $\fp\in T$, we explain how to extend the
computation of $\pi_{\fC}$ on $\Gamma_{T\setminus \{\fp\}}$ to the map
\[
    \rho_{\fC}:\Gamma_T \longrightarrow \wtU_\fp^k.
\]

\subsection{The permutation representation}
For the purposes of computing $\pi_\fC$, we realize $G[\fC]$ as a
permutation group. Let
\[
    \overline{\fo}_{\fC}=\fo_K/\fC,
    \qquad
    \overline{\fO}_{\fC}=\fO/\fC\fO.
\]
Reducing the structure constants $c_{ij\ell}$ modulo $\fC$ gives an
explicit $\overline{\fo}_{\fC}$-algebra structure on
$\overline{\fO}_{\fC}$. We then enumerate the finite unit group
$\overline{\fO}_{\fC}^{\,\times}$ and its scalar subgroup
$\overline{\fo}_{\fC}^{\,\times}$, and form the quotient set
\[
    X_\fC =
    \overline{\fO}_{\fC}^{\,\times}/\overline{\fo}_{\fC}^{\,\times}.
\]
Left multiplication on $X_\fC$ gives a faithful permutation
representation
\[
    \sigma_\fC:G[\fC]\hookrightarrow \Sym(X_\fC).
\]
Given representatives $\alpha \in \fO[1/T]$ for the elements of
$\Gamma_T$, we reduce their coordinates in the basis
$\omega_1,\ldots,\omega_4$ modulo $\fC$, using that the primes in $T$
are invertible modulo $\fC$. This gives
\[
    \overline{\Gamma}_{\fC,T}
    =
    \sigma_\fC(\pi_\fC(\Gamma_T)).
\]
Since $\sigma_\fC$ is injective, the map $\pi_\fC$ is injective on
$\Gamma_T$ if and only if
\[
    \#\overline{\Gamma}_{\fC,T}=\#\Gamma_T.
\]
Finally, we check in Magma whether the finite permutation group
$G[\fC]$ contains a complement of $\overline{\Gamma}_{\fC,T}$.

In practice, it suffices to construct $\sigma_{\fC}$ when $\fC$ is a
prime power. The general case then follows by taking direct products of
the prime-power computations via Sunzi's remainder theorem; see
\eqref{eq:sunzi-factorization}.

\begin{rmk}
    Since $\fo_K/\fp$ is a finite field, Magma can test the
    invertibility of $\overline\alpha \in \overline{\fO}_\fp$
    directly. For $\fC=\fp^k$ with $k\geq1$, the ideal
    $\fp\overline{\fO}_{\fC}$ is nilpotent. Thus,
    $\overline\alpha \in \overline{\fO}_{\fC}$ is a unit if and only
    if the reduction of $\overline\alpha$ in $\overline{\fO}_\fp$ is a
    unit. Therefore, computing $\overline{\fO}_{\fC}^\times$ reduces to
    lifting elements of $\overline{\fO}_{\fp}^\times$.
\end{rmk}

\begin{rmk}
    Magma has a routine \texttt{UnitGroup} which computes
    representatives for
    \[
        \Gamma_{\emptyset}=\fO^\times/\fo_K^\times,
    \]
    and a routine \texttt{Normalizer} which computes representatives for
    \[
        N_{\emptyset}=N_{B^\times}(\fO)/K^\times.
    \]
    Since $\Gamma_{\emptyset} \subseteq \Gamma_T\subseteq N_{\emptyset}$
    for every $T\subseteq \Ram_f(B)$, these routines help us compute
    $\Gamma_T$. This computation is easiest when $\fO[1/T]$ is Eichler,
    given the explicit description of the quotient $N_T/\Gamma_T$ in
    Proposition~\ref{prop:normalizer-unit-exact-sequence}.
\end{rmk}

\subsection{Extending the representation}

Now let
\[
    T \subseteq \Ram_f(B), \quad
    \fp \in \Ram_f(B) \setminus T, \quad\text{and}\quad
    S = T \cup \{\fp\}.
\]
We further suppose that $\fO[1/T]$ is Eichler with class number one.
Let $\fC = \fp^k$ and identify $U^{k}_\fp$ with $G[\fC]$ via
$\bar\eta$. Recall from the paragraph after Proposition \ref{prop:normalizer-unit-exact-sequence}
that $\rho_\fC$ is injective on $\Gamma_S$ if and only if it is injective on $\Gamma_T$;
see also Proposition \ref{prop:ramified-complement-lift}.
With the methods of the previous subsection, we compute the
permutation representation
\[
    \sigma_\fC: G[\fC] \hookrightarrow \Sym(X_\fC)
\]
and the image $\overline{\Gamma}_{\fC,T}$. We compute the $C_2$-extension
$\wtU_\fp^{k}$ of $U_\fp^k$ and the map $\rho_\fC$ on $\Gamma_S$ as follows.

Let $\fP$ be the unique two-sided $\fO$-ideal of reduced norm $\fp$,
and choose a generator
\[
    \fP = \alpha \cdot \fO[1/T].
\]
Then $\alpha$ represents the nontrivial coset in
$\Gamma_S/\Gamma_T \cong C_2$ (see
Proposition \ref{prop:normalizer-unit-exact-sequence}), and $\alpha^2 \in \Gamma_T$.
Set
\[
    \gamma_\alpha
    =
    \sigma_\fC(\pi_\fC(\alpha^2))
    \in
    \overline{\Gamma}_{\fC,T}.
\]
Even though the short exact sequence \eqref{eq:ramified-projective-unit-extension} splits,
$\gamma_\alpha$ may not be trivial, so we cannot simply write
$\wtU_\fp^k$ as a semidirect product using the chosen lift $\alpha$.
However, $\alpha$ acts on $G[\fC]$ by conjugation. Let
\[
    \psi_\alpha:G[\fC]\to G[\fC],
    \qquad
    \psi_\alpha(g)=\alpha g\alpha^{-1},
\]
denote the induced automorphism. Since $\overline{\Gamma}_{\fC,T}$ is
stable under $\psi_\alpha$, after computing presentations for the
permutation groups $\overline{\Gamma}_{\fC,T}$ and $G[\fC]$, we compute presentations
\[
    \overline{\Gamma}_{\fC,S} = \rho_\fC(\Gamma_S)
    \cong
    \lag
        \overline{\Gamma}_{\fC,T}, \lambda_\alpha
        \ \middle|\
        \lambda_\alpha^2 = \gamma_\alpha,\
        \lambda_\alpha g \lambda_\alpha^{-1}=\psi_\alpha(g)
        \text{ for } g\in \overline{\Gamma}_{\fC,T}
    \rag,
\]
and
\[
    \wtU_\fp^k
    \cong
    \lag
        G[\fC], \lambda_\alpha
        \ \middle|\
        \lambda_\alpha^2 = \gamma_\alpha,\
        \lambda_\alpha g \lambda_\alpha^{-1}=\psi_\alpha(g)
        \text{ for } g\in G[\fC]
    \rag.
\]
Having computed these finite groups, we then look for complements of
$\overline{\Gamma}_{\fC,S}$ in $\wtU_\fp^k$.

\begin{rmk}
    Given a finite group $G$, Magma lets one easily switch between a
    permutation representation of $G$ and a finite presentation of $G$.
    Thus, it is reasonable to compute $\overline{\Gamma}_{\fC,S}$ and
    $\wtU_\fp^k$ in the above fashion.
\end{rmk}

\subsection{Classification of minimal complementary ideals}

For the remainder of this section, we assume that $\fO \subseteq B$ is maximal
and class number one. In \cite{KirschmerVoight10},
Kirschmer and Voight classify the totally definite quaternion algebras
over number fields with a maximal order of class number one. There are
49 such algebras, and in each case $[K:\bQ] \leq 5$.

We would like to classify all complementary triples associated to a
class number one maximal order. Although such a classification may be
attainable (see Section \ref{sec:irreducible-complementary-triples}), we do not
complete it in this paper. Instead, we restrict our attention to a more
tractable problem.

We call a complementary triple $(T,\fC,H)$ of $\fO$ \emph{minimal} if,
for every proper divisor $\fa \mid \fC$, the corresponding reduction map
is not injective on $\Gamma_T$; here this map is $\pi_\fa$ in the
uninverted case and $\rho_\fa$ in the inverted case. In this section, we
give a method for classifying all minimal complementary triples of
$\fO$. Our classification relies on a classical result of Dickson
\cite{Dickson01} classifying all maximal subgroups of
$\PGL_2(\bF_q)$ for prime powers $q$. See \cite{King05} for a modern
exposition of Dickson's theorem. In particular, we use
\cite[Corollary 2.3]{King05} to derive the following lemma.

\begin{lem} \label{lem:dickson-index-bound}
    Let $q$ be a prime power and let $H \leq \PGL_2(\bF_q)$. Then
    \[
        [\PGL_2(\bF_q):H] \geq q+1
    \]
    unless one of the following holds:
    \begin{enumerate}
        \item $H = \PSL_2(\bF_q)$ of index $2$;
        \item $q = 5$ and $H$ is of index $5$ with $H \cong S_4$
        under the isomorphism $\PGL_2(\bF_5) \cong S_5$;
        \item $q = 3$ and $H$ is of index $3$ with $H \cong D_4$
        under the isomorphism $\PGL_2(\bF_3) \cong S_4$.
    \end{enumerate}
\end{lem}

We now confirm the following lemma with Magma.

\begin{lem} \label{lem:unit-representatives-norm-one}
    Every element of $\Gamma_\emptyset$ has a lift to $\fO^\times$ with
    reduced norm $1$.
\end{lem}

The next lemma completely determines when $\fO$ has infinitely many
minimal complementary triples.

\begin{lem} \label{lem:cyclic-unit-trivial-complementary-triples}
    If $\Gamma_\emptyset = 1$, then $(\emptyset,\fC,G[\fC])$ is a
    complementary triple of $\fO$ for every ideal $\fC$. Furthermore, if
    $\#T = 1$ and $\alpha$ is a generator of
    $\Gamma_T \cong C_2$, then the following statements hold.
    \begin{enumerate}[label=(\roman*)]
        \item \label{case:cyclic-same-ramified-prime} If $\fp$ ramifies and
        $T = \{\fp\}$, then for every $k \geq 1$,
        \[
            (T,\fp^k,U_\fp^k)
        \]
        is a complementary triple of $\fO$.

        \item \label{case:cyclic-minimal-prime-power} If $(T,\fC,H)$ is a
        minimal complementary triple of $\fO$, then $\fC = \fp^k$ for
        some prime ideal $\fp$. Furthermore, $k = 1$ if $\fp$ is split
        or $\fp \in T$.

        \item \label{case:cyclic-distinct-ramified-primes} If $\fp$ and $\fq$
        are distinct ramified primes and $T = \{\fq\}$, then
        $(T,\fp^k,H)$ is a minimal complementary triple if and only if
        \[
            K = \bQ(\sqrt{2}), \qquad
            \fD(B) = \fp\fq, \qquad
            k=1,
        \]
        where $\fp$ and $\fq$ are the unique primes of norms $25$ and
        $2$, respectively.

        \item \label{case:cyclic-odd-split-prime} Let $\fp$ be a split prime
        with odd norm, and let $H \leq G[\fp]$. Then $(T,\fp,H)$ is a
        complementary triple of $\fO$ if and only if $\Nrd(\alpha)$ is
        not a square modulo $\fp$, and
        \[
            H = \PSL_2(\fo_K/\fp).
        \]

        \item \label{case:cyclic-even-split-prime} Let $\fp$ be a split prime
        with even norm, and let $H \leq G[\fp]$. Then $(T,\fp,H)$ is a
        complementary triple of $\fO$ if and only if $N\fp = 2$, and,
        under the identification
        \[
            G[\fp] \cong \PGL_2(\bF_2) \cong S_3,
        \]
        the subgroup $H$ is the unique subgroup isomorphic to $A_3$.
    \end{enumerate}
\end{lem}

\begin{proof}
    The first statement is immediate. Suppose that $\fp$ ramifies and
    $k \geq 1$. Since $\Gamma_\emptyset = 1$, applying the fact that
    $(\emptyset,\fp^k,G[\fp^k])$ is a complementary triple with
    Proposition \ref{prop:ramified-complement-lift} proves
    \ref{case:cyclic-same-ramified-prime}.

    Now let $(T,\fC,H)$ be a minimal complementary triple of $\fO$, and
    write
    \[
        \fC = \prod_{i=1}^r \fp_i^{k_i}.
    \]
    We first show that $r=1$. Suppose, for contradiction, that
    $r \geq 2$. Factor $G[\fC]$ as in \eqref{eq:sunzi-factorization}. Since
    $\#\Gamma_T = 2$, the subgroup $H$ has index two in $G[\fC]$.
    Hence, $H$ is the kernel of a nontrivial map
    \[
        \chi : G[\fC] \to C_2
    \]
    such that
    \[
        \chi(\pi_\fC(\alpha)) \neq 1.
    \]
    Since $\chi$ is determined by its restrictions to the factors
    $G[\fp_i^{k_i}]$, there exists some $i$ for which the restriction of
    $\chi$ to $G[\fp_i^{k_i}]$ is nontrivial. Therefore,
    \[
        \pi_{\fp_i^{k_i}}(\alpha) \neq 1,
    \]
    which contradicts the minimality of $\fC$. Thus, $r=1$, proving the
    first claim of \ref{case:cyclic-minimal-prime-power}.

    We now verify the second claim. If $\fp \in T$, then
    \ref{case:cyclic-same-ramified-prime} applied with $k=1$ shows that the
    reduction map at $\fp$ is already injective on $\Gamma_T$, so
    minimality forces $k=1$. If $\fp$ splits, then we check in Magma
    that $\pi_\fp$ must be injective on $\Gamma_T$, and again minimality
    forces $k=1$.

    Given the set-up of \ref{case:cyclic-distinct-ramified-primes}, we verify
    in Magma that, among the finitely many algebras in the
    Kirschmer--Voight classification, the given algebra $B$ is the only
    one with $\#\Ram_f(B)>1$ and $\Gamma_\emptyset=1$. Moreover, the
    triple described is the only minimal triple of the form
    $(\{\fq\},\fp^k,H)$.

    Now let $\fp$ be a split prime. Then $\pi_\fp$ is injective on
    $\Gamma_T$ by \ref{case:cyclic-minimal-prime-power}. Since $\fp$ is split,
    we identify
    \[
        G[\fC] \cong \PGL_2(\fo_K/\fp^k).
    \]
    Let $(T,\fp,H)$ be a complementary triple of $\fO$.

    Suppose that $N\fp$ is odd. If $\Nrd(\alpha)$ is a square modulo
    $\fp^k$, then the image of $\Gamma_T$ lies in
    $\PSL_2(\fo_K/\fp^k)$. Since $H$ is an index two complement of
    $\Gamma_T$, the subgroup
    \[
        H \cap \PSL_2(\fo_K/\fp)
    \]
    would then be an index two subgroup of $\PSL_2(\fo_K/\fp)$, but
    $\PSL_2(\fo_K/\fp)$ has no subgroup of index two. Hence,
    $\Nrd(\alpha)$ cannot be a square modulo $\fp$.

    Conversely, suppose that $\Nrd(\alpha)$ is not a square modulo
    $\fp$. Then the nontrivial element of $\pi_\fC(\Gamma_T)$ lies
    outside $\PSL_2(\fo_K/\fp)$. Since $\PSL_2(\fo_K/\fp)$ is the
    unique index two subgroup of $\PGL_2(\fo_K/\fp)$, it follows that
    $H = \PSL_2(\fo_K/\fp)$. This proves
    \ref{case:cyclic-odd-split-prime}.

    If $\fp$ has even norm, then
    $\PGL_2(\fo_K/\fp) \cong \PSL_2(\fo_K/\fp)$ contains an index two
    subgroup if and only if $N\fp = 2$, in which case
    \[
        G[\fC] \cong \PGL_2(\bF_2) \cong S_3.
    \]
    The unique index two subgroup of $S_3$ is $A_3$. Since the
    nontrivial image of $\Gamma_T$ has order two, it lies outside
    $A_3$, and thus $A_3$ complements $\pi_\fC(\Gamma_T)$. This proves
    \ref{case:cyclic-even-split-prime}.
\end{proof}

A complementary triple $(T,\fC,H)$ of $\fO$ is called \emph{trivial} if
$\Gamma_\emptyset = 1$, $\#T \in \{0,1\}$, and $(T,\fC,H)$ is not the
exceptional triple described in
Lemma~\ref{lem:cyclic-unit-trivial-complementary-triples}\ref{case:cyclic-distinct-ramified-primes}. Thus, all
triples described in Lemma~\ref{lem:cyclic-unit-trivial-complementary-triples} except the one in
Lemma~\ref{lem:cyclic-unit-trivial-complementary-triples}\ref{case:cyclic-distinct-ramified-primes} are trivial.
We exclude this exceptional triple from the trivial class because it is
not forced by some formal argument as in the other cases.

We record one more lemma before our main classification theorem for
nontrivial, minimal complementary triples.

\begin{lem} \label{lem:c2-unit-trivial-complementary-triples}
    If $\Gamma_\emptyset \cong C_2$ and $(\emptyset, \fC, H)$ is a
    minimal triple of $\fO$, then $\fC = \fp^k$ and either $\fp$
    ramifies or $\fp$ splits with $N\fp = 2$ and $k = 1$.
\end{lem}

\begin{proof}
    A similar argument to
    Lemma \ref{lem:cyclic-unit-trivial-complementary-triples}\ref{case:cyclic-minimal-prime-power} shows that
    $\fC = \fp^k$. Suppose that $\fp$ splits. We check in Magma that
    $\pi_\fp$ must be injective on $\Gamma_\emptyset$. Under the isomorphism
    \[
        G[\fp] \cong \PGL_2(\fo_K/\fp),
    \]
    we see from Lemma \ref{lem:unit-representatives-norm-one} that
    $\pi_\fC(\Gamma_\emptyset) \subseteq \PSL_2(\fo_K/\fp)$, and
    \[
        H \cap \PSL_2(\fo_K/\fp)
    \]
    is an index two subgroup of $\PSL_2(\fo_K/\fp)$. Such a subgroup
    exists if and only if $N\fp = 2$.
\end{proof}

\begin{thm} \label{thm:minimal-complementary-triple-classification}
    Let $(T,\fC,H)$ be a nontrivial, minimal complementary triple of
    $\fO$. If $\fp \mid \fC$ is prime, then $\fp$ is ramified or
    \[
        N\fp \leq \#\Gamma_T.
    \]
    Equality occurs if and only if one of the following holds:
    \begin{enumerate}
        \item $\Gamma_T \cong C_2$ and $\fC = \fp$ with $N\fp = 2$;
        \item $\Gamma_T \cong C_5$ and $\fC = \fp$ with $N\fp = 5$.
    \end{enumerate}
    In particular, there are only finitely many minimal nontrivial
    complementary triples of $\fO$.
\end{thm}

\begin{proof}
    Let $\fp \mid \fC$ be a split prime, and suppose for contradiction
    that $N\fp>\#\Gamma_T$. Using Magma, we write representatives for
    $\Gamma_T$ in terms of the basis in \eqref{eq:order-basis} and check
    that $\pi_\fp$ is injective on $\Gamma_T$. By the minimality of
    $\fC$, it follows that $\fC=\fp$. The result now follows from
    Lemma \ref{lem:dickson-index-bound} and Lemma \ref{lem:c2-unit-trivial-complementary-triples}.
\end{proof}

\begin{ex}
    It follows from Theorem \ref{thm:minimal-complementary-triple-classification} that, if $\fp$
    splits and $N\fp > 5$ is even, then $N\fp < \#\Gamma_T$. This bound
    is optimal. For instance, the maximal order $\fO$ of
    $B = \quat{\bQ(\sqrt{5})}{-1}{-1}$ has
    $\Gamma_{\fO,\emptyset} \cong A_5$, and the primes of norm $59$ are
    complementary primes of $\fO$ at $\emptyset$.
\end{ex}

\begin{ex}
    We give another curious edge case. Let
    \[
        K = \bQ(\zeta+\zeta^{-1}), \quad \zeta = e^{2\pi i/7}.
    \]
    Then there is a unique prime $\fp_7 \subseteq \fo_K$ of norm $7$.
    There is a totally definite quaternion algebra $B/K$ of class number
    one with discriminant $\fD(B) = \fp_7$. Let $\fO \subseteq B$ be a
    maximal order. Setting $T = \{\fp_7\}$, we have
    \[
        \Gamma_\emptyset \cong D_7, \quad \Gamma_T = D_{14}.
    \]
    It will follow from our classification of minimal triples in Appendix \ref{app:complementary-triple-tables}
    that the only minimal complementary triple of $\fO$ is $(T,\fp_{27},H)$,
    where $\fp_{27}$ is the prime of norm $27$.
\end{ex}

Let $\cM(B)$ be the set of nontrivial, minimal complementary triples of
the maximal order $\fO \subseteq B$, up to equivalence.
In Appendix \ref{app:complementary-triple-tables}, we use
Theorem \ref{thm:minimal-complementary-triple-classification} to enumerate $\cM(B)$ for each totally
definite quaternion algebra $B/K$ of class number one enumerated by
Kirschmer and Voight. Our computations yield the following result.

\begin{thm} \label{thm:nonempty-minimal-triples}
    Let $K$ be a number field, and let $B/K$ be a totally definite
    quaternion algebra of class number one with
    $\Gamma_\emptyset \neq 1$ or $\#\Ram_f(B) > 1$.
    Then $\cM(B)$ is non-empty.
\end{thm}

\begin{cor} \label{cor:existence-complementary-triples}
    Let $K$ be a number field, and let $B/K$ be a totally definite
    quaternion algebra of class number one with maximal order $\fO$.
    Then there exists a complementary triple of $\fO$.
\end{cor}

\begin{cor}
    Theorems \ref{thm:cayley-ramanujan-graphs}, \ref{thm:cayley-cubical-complexes}, and \ref{thm:same-valency-products} follow.
\end{cor}

\begin{rmk} \label{rem:complementary-triples-necessary}
    Suppose that $B/K$ has class number one and maximal order $\fO$.
    Any order $\fL \subseteq \fO$ with class number one and
    $\Gamma_{\fL,\emptyset} = 1$ immediately yields simply transitive
    actions $\Gamma_{\fL,S} \leq \Gamma_{\fO,S}$ on $\cT_S$ for
    infinitely many choices of $S$. One might hope that, for every
    quaternion algebra $B/K$ of class number one, there exists such an
    order of class number one. If so, the search for complementary
    triples would seem unnecessarily arduous.

    However, it follows from \cite{KirschmerLorch16} that there are
    class number one quaternion algebras $B/K$ for which the only class
    number one order in $B$ is the maximal order $\fO$. Thus,
    complementary triples are necessary in order to construct subgroups
    of $\Gamma_{\fO,S}$ that act simply transitively on
    $\cG_S(\fO)$.
\end{rmk}

We give one more interesting outcome of our classification.

\begin{prop}
    Let $K$ be a number field, let $B/K$ be a totally definite
    quaternion algebra of class number one, and let
    $(T,\fC,H) \in \cM(B)$. Then $\fC$ is a prime power except when
    \[
        K = \bQ(\zeta + \zeta^{-1}), \qquad
        \zeta = e^{2\pi i/9}, \qquad
        \fD(B) = \fp_3,
    \]
    where $\fp_3$ is the unique prime above $3$. In this exceptional
    case, the only element of $\cM(B)$ is represented by a triple
    $(T,\fC,H)$ with
    \[
        T = \{\fp_3\}
        \qquad\text{and}\qquad
        N(\fC) = 24.
    \]
\end{prop}

\section{Torsion-free groups} \label{sec:torsion-free-groups}

In this section, we study when the group $\Gamma_{\fO,S}(\fC,H)$
associated to a complementary triple is torsion-free.
We begin with a general lemma about torsion in simply transitive lattices.

\begin{lem} \label{lem:simply-transitive-torsion-order-two}
    If $\Gamma$ acts simply transitively on $\cG_S(\fO)$, then every nontrivial torsion element of $\Gamma$ has order $2$.
\end{lem}

\begin{proof}
    Let $\gamma \in \Gamma$ be a nontrivial torsion element. By the
    Bruhat--Tits fixed point lemma, $\gamma$ fixes a cell $\sigma$ of
    $\cT_{S_0}$; see \cite[II.2.8(1)]{BridsonHaefliger99}. Since
    $\Gamma$ acts simply transitively on vertices, $\gamma$ cannot fix a
    vertex, so $\sigma$ has positive dimension.

    Suppose $\sigma$ is an $n$-cube, and fix an identification with
    $[0,1]^n$. The stabilizer of $\sigma$ embeds into the automorphism
    group of the $n$-cube, which is isomorphic to
    \[
        C_2^n\rtimes S_n.
    \]
    The symmetric group $S_n$ permutes the coordinates of $[0,1]^n$.
    However, each coordinate comes from a unique tree factor
    $\cT(K_\fp)$, and elements of
    $\Gamma \leq \prod_{\fp\in S_0}\PGL_2(K_\fp)$ do not permute the
    factors. Therefore, the stabilizer of $\sigma$ embeds in $C_2^n$,
    so $\gamma$ must have order $2$.
\end{proof}

\subsection{Torsion obstructions}

Let $K$ be a number field, let $B/K$ be a totally definite quaternion
algebra, and let $\fO \subseteq B$ be a not necessarily maximal order.
Fix a complementary triple $(T,\fC,H)$ of $\fO$, and let $S$ be
$\fC$-eligible. Put
\[
    \Gamma_{\fO,S}(\fC,H) = \pi_S^{-1}(H)
\]
as before, and assume that $\Gamma_{\fO,S}(\fC,H)$ acts transitively on
$\cG_S(\fO)$. Then Proposition~\ref{prop:complementary-triple-simple-transitive}
implies that $\Gamma_{\fO,S}(\fC,H)$ acts simply transitively on
$\cG_S(\fO)$, so Lemma \ref{lem:simply-transitive-torsion-order-two} shows that the only possible
nontrivial torsion in $\Gamma_{\fO,S}(\fC,H)$ is $2$-torsion.

Recall that we have a basis $\omega_1,\ldots,\omega_4$ of $\fO$ and a
quadratic form
\[
    Q_\fO(X_1,X_2,X_3,X_4)
    =
    \sum_i \Nrd(\omega_i)X_i^2
    +
    \sum_{j<k} \trd(\omega_j\overline{\omega_k}) X_jX_k
\]
with coefficients in $\fo_K$ such that
\[
    \Nrd(\alpha) = Q_\fO(X_1,X_2,X_3,X_4),
    \quad
    \alpha = \sum_{i=1}^4 X_i\omega_i.
\]
Let $\fa \subseteq \fo_K$ be an ideal and $\Delta_\fa = \Delta(\fo_K/\fa)$. The reduced norm induces a well-defined map
\[
    \Nrd_\fa : G[\fa] \to \Delta_\fa
\]
which fits into the following commutative diagram.
\[\begin{tikzcd}
	{\Gamma_{\fO,S}} && {\Delta(\fo_K[1/S])} \\
	{G[\fa]} && {\Delta_\fa}
	\arrow["\Nrd", from=1-1, to=1-3]
	\arrow["{\pi_\fa}"', from=1-1, to=2-1]
	\arrow["{\pi_\fa}", from=1-3, to=2-3]
	\arrow["{\Nrd_\fa}"', from=2-1, to=2-3]
\end{tikzcd}\]
Note that if $\fa = \fp^k$ for $\fp$ ramifying in $B$ and $\fO_\fp$ is maximal, then we have recovered the form described in Section \ref{sec:nonsplit-local-case}.

Observe that the trace form also reduces down to $\fO/\fa\fO$. Given
$g \in G[\fa]$, we say that $g$ has \emph{trace-zero}, denoted
\[
    \trd_\fa(g) = 0,
\]
if some lift $\overline\alpha \in (\fO/\fa \fO)^\times$ of $g$ has trace
$0$. This condition is independent of the choice of lift.

For $Z\subseteq G[\fa]$, define the \emph{torsion obstruction set at
$Z$} to be
\[
    \Theta_{\fa}(Z)
    =
    \left\{
        \Nrd_{\fa}(g)
        :
        g\in Z,\ \trd_{\fa}(g)=0
    \right\}
    \subseteq \Delta_\fa.
\]
When $Z=G[\fa]$, we write
\[
    \Theta_{\fa}=\Theta_{\fa}(G[\fa]).
\]
Given our computations in the previous section, the set $\Theta_\fa(Z)$
is straightforward to compute. We now show how the torsion obstruction
set gives a criterion for torsion-freeness.

\begin{prop}[Torsion obstruction criterion] \label{prop:torsion-obstruction-criterion}
    Let $\fa\subseteq \fo_K$ be an ideal coprime to every finite prime in $S$, and set
    \[
        \Gamma=\Gamma_{\fO,S}(\fC,H),
        \qquad
        Z_\fa=\pi_\fa(\Gamma)\subseteq G[\fa],
        \qquad
        \Delta_{\fa,S}
        =
        \pi_\fa\bigl(\Delta^+(\fo_K[1/S])\bigr)
        \leq \Delta_\fa.
    \]
    If $\gamma\in \Gamma$ is a nontrivial element of order $2$, then
    \[
        \Nrd_\fa\bigl(\pi_\fa(\gamma)\bigr)
        \in
        \Delta_{\fa,S}\cap \Theta_\fa(Z_\fa).
    \]
    Consequently, if
    \begin{equation} \label{eq:torsion-obstruction-emptiness}
        \Delta_{\fa,S}\cap \Theta_\fa(Z_\fa)=\emptyset
    \end{equation}
    then $\Gamma_{\fO,S}(\fC,H)$ is torsion-free.
\end{prop}

\begin{proof}
    Let $\gamma\in\Gamma$ be a nontrivial element of order $2$, and
    choose a representative $\alpha\in \fO[1/S]^\times$. Since
    $\gamma$ has order $2$, we have
    $\alpha^2\in \fo_K[1/S]^\times$, while $\alpha$ is not scalar. The
    well-known identity
    \[
        \alpha^2-\trd(\alpha)\alpha+\Nrd(\alpha)=0
    \]
    forces $\trd(\alpha)=0$. After reducing modulo $\fa$, this shows
    that $\pi_\fa(\gamma)\in Z_\fa$ is represented by a trace-zero
    element. Hence,
    \[
        \Nrd_\fa\bigl(\pi_\fa(\gamma)\bigr)
        \in \Theta_\fa(Z_\fa).
    \]

    On the other hand, since $\alpha\in\fO[1/S]^\times$, its reduced
    norm lies in $\fo_K[1/S]^\times$. Moreover, $B$ is totally
    definite, so $\Nrd(\alpha)$ is totally positive. Thus,
    $\Nrd(\alpha)$ represents an element of
    $\Delta^+(\fo_K[1/S])$, and reducing modulo $\fa$ gives
    \[
        \Nrd_\fa\bigl(\pi_\fa(\gamma)\bigr)
        \in \Delta_{\fa,S}.
    \]
    By Lemma \ref{lem:simply-transitive-torsion-order-two}, every nontrivial torsion element of
    $\Gamma$ has order $2$. Hence, $\Gamma$ is torsion-free if
    \eqref{eq:torsion-obstruction-emptiness} holds.
\end{proof}

\begin{cor} \label{cor:empty-obstruction-set-torsion-free}
    If $\Theta_\fa(Z_\fa)=\emptyset$, then
    $\Gamma_{\fO,S}(\fC,H)$ is torsion-free.
\end{cor}

\begin{rmk}
    The subgroup $\Delta_{\fa,S}$ is easy to compute in practice. If
    each $\fp\in S$ has a chosen totally positive generator $a_\fp$,
    then
    \[
        \Delta_{\fa,S}
        =
        \left\langle
            \pi_\fa(u),\ \pi_\fa(a_\fp)
            :
            [u]\in \Delta^+(\fo_K),\ \fp\in S
        \right\rangle.
    \]
\end{rmk}

The torsion obstruction set $\Theta_\fa(Z_\fa)$ is only useful when it
does not contain $1$. Otherwise, its intersection with $\Delta_{\fa,S}$
is automatically non-empty. The next proposition shows that, for
maximal orders, useful torsion obstructions can only occur at primes
dividing the discriminant or the complementary ideal.

\begin{prop} \label{prop:torsion-obstruction-prime-support}
    Let $S$ be $\fC$-eligible, and let $\fa\subseteq \fo_K$ be an ideal
    coprime to $\fD(B)$, $\fC$, and every finite prime in $S$. Assume
    that $\fO_\fq$ is maximal for every $\fq\mid \fa$. Then
    \[
        \Theta_{\fa}= \Delta_\fa \quad \text{and} \quad
        1 \in \Theta_\fa(Z_\fa).
    \]
    In particular, $\Delta_{\fa,S} \cap \Theta_\fa(Z_\fa)$ is non-empty.
\end{prop}

\begin{proof}
    Since $\fO_\fq$ is maximal and $B$ splits at every $\fq\mid \fa$, we have an isomorphism $G[\fa] \cong \PGL_2(\fo_K/\fa)$
    sending the reduced norm to the determinant and trace-zero elements to trace-zero elements.
    Let $c\in \Delta(\fo_K/\fa)$. The matrix
    \[
        \begin{pmatrix}
            0 & -c \\
            1 & 0
        \end{pmatrix}
    \]
    has trace-zero and determinant $c$. Hence, every class in
    $\Delta_\fa$ occurs in $\Theta_{\fa}$, proving the first claim.

    It remains to show that the trivial class occurs in $\Theta_\fa(Z_\fa)$. Consider
    \[
        J =
        \begin{pmatrix}
            0 & -1 \\
            1 & 0
        \end{pmatrix}
        \in \PSL_2(\fo_K/\fa).
    \]
    This matrix has trace-zero and determinant $1$. Since $\fa$ and
    $\fC$ are coprime, strong approximation for $\fO[1/S]^1$
    \cite[Main Theorem 28.5.3]{Voight21} tells us that we may choose
    $\alpha \in \fO[1/S]^1$ such that
    \[
        \pi_\fa(\alpha) = J, \quad \pi_\fC(\alpha) = 1.
    \]
    Therefore, $\alpha$ represents an element in
    $\Gamma_{\fO,S}(\fC,H)$ such that
    \[
        \trd_\fa(\pi_\fa(\alpha))=0
        \qquad\text{and}\qquad
        \Nrd_\fa(\pi_\fa(\alpha))=1,
    \]
    proving the second claim.
\end{proof}

We now study the obstruction sets occurring at the complementary prime
$\fC$ and the discriminant $\fD(B)$. First, we will need the following
lemma.

\begin{lem}
    Let $\fp \subseteq \fo_K$ be a prime ideal, and set
    \[
        m = 2v_{\fp}(2)+1.
    \]
    Then, for every $k \geq m$, the reduction map
    $\fo_K/\fp^k \to \fo_K/\fp^m$ induces an isomorphism
    \[
        \Delta_{\fp^k} \cong \Delta_{\fp^m}.
    \]
\end{lem}

\begin{proof}
    Surjectivity is immediate. It remains to prove injectivity. Let
    $\bar a,\bar b \in (\fo_K/\fp^k)^\times$ have the same image in
    $\Delta_{\fp^m}$. Choose lifts $a,b\in \fo_{K,\fp}^\times$. Then
    there exists $x\in \fo_{K,\fp}^\times$ such that
    \[
        ab^{-1} \equiv x^2 \pmod{\fp^m}.
    \]
    Set $c=ab^{-1}$ and consider
    \[
        f(X)=X^2-c.
    \]
    Then
    \[
        v_{\fp}(f(x)) \geq m = 2v_{\fp}(2)+1.
    \]
    Since $x\in \fo_{K,\fp}^\times$, we have
    \[
        v_{\fp}(f'(x))
        =
        v_{\fp}(2x)
        =
        v_{\fp}(2).
    \]
    Hence,
    \[
        v_{\fp}(f(x))
        >
        2v_{\fp}(f'(x)).
    \]
    By Newton's lemma \cite[Theorem 7.32]{MilneANT}, $f$ has a root
    $y\in \fo_{K,\fp}$. Reducing modulo $\fp^k$, we see that
    $\bar a\bar b^{-1} = \bar y^2$ is a square in
    $(\fo_K/\fp^k)^\times$. Therefore $\bar a$ and $\bar b$ represent
    the same class in $\Delta_{\fp^k}$.
\end{proof}

Recall that by the Sunzi remainder theorem, we have
\[
    \Delta_\fa \cong \prod_{i=1}^r \Delta_{\fp_i^{e_i}},
    \qquad
    \fa = \fp_1^{e_1} \cdots \fp_r^{e_r}.
\]
Let
\[
    \fe(\fa) = \prod_{\fp \mid \fa} \fp^{m(\fa,\fp)}, \qquad
    m(\fa,\fp) =
    \max\left\{2v_\fp(2)+1,v_\fp(\fa)\right\}.
\]
It follows from the previous lemma that, if $\fa\mid \fb$ and $\fa$
and $\fb$ have the same prime divisors, then
$\Delta_{\fe(\fa)} \cong \Delta_{\fe(\fb)}$.
We now obtain a corollary of Proposition \ref{prop:torsion-obstruction-prime-support}.

\begin{cor} \label{cor:ramified-prime-obstruction-suffices}
    Let $S$ be $\fC$-eligible, let $(T,\fC,H)$ be a complementary
    triple, and let
    \[
        \widetilde H = \pi_{\fe(\fC),\fC}^{-1}(H),
    \]
    where $\pi_{\fe(\fC),\fC}:G[\fe(\fC)]\to G[\fC]$ is the natural
    reduction map. If
    $1 \in \Theta_{\fe(\fC)}(\widetilde H)$ and
    $1 \in \Theta_{\fe(\fD(B))}$, then
    $1 \in \Theta_\fa(Z_\fa)$ for all ideals $\fa \subseteq \fo_K$
    coprime to every finite prime in $S$.
\end{cor}

In Appendix \ref{app:complementary-triple-tables}, we determine whether
$1 \in \Theta_{\fe(\fD(B))}$ and whether
$1 \in \Theta_{\fe(\fC)}(\widetilde H)$ for each class number one
algebra $B$ in the Kirschmer--Voight classification and for each triple
$(T,\fC,H) \in \cM(B)$.

\begin{ex}
    There exist complementary triples $(T,\fC,H)$ with $\fC$-eligible
    $S$ such that
    \[
        \Theta_\fa(Z_\fa) \cap \Delta_{\fa,S} \neq \emptyset
    \]
    for all $\fa$, but $\Gamma_{\fO,S}(\fC,H)$ is torsion-free. We saw
    one such case in Example \ref{ex:qsqrt2-torsion-free}.

    Indeed, we showed via geometric methods that the lattice
    $\Gamma = \Gamma_{\fO,S}(\fC,H)$ given in this example is
    torsion-free. However, in
    Appendix \ref{app:complementary-triple-tables}, we see that
    $1 \in \Theta_{\fe(\fp_2)}(\widetilde H)$, so the intersection with
    $\Delta_{\fe(\fp_2)}$ is non-empty. Furthermore, in Magma we find
    \[
        \Theta_{\fp_9}=\left\{[1+\sqrt2]\right\}.
    \]
    Identifying $\fo_K/\fp_9$ with $\bF_3[t]/(t^2-2)$, the square units
    are given by $\left\{1,2,t,2t\right\}$. Consequently, the
    fundamental unit $1+\sqrt2$ of $\fo_K$ is not a square modulo
    $\fp_9$. Therefore, given an ideal $\fp \in S$, we can pick a
    generator $a$ of $\fp$ with $\pi_{\fp_9}(a) = [1 + \sqrt{2}]$, so
    the intersection with $\Delta_{\fa,S}$ is non-empty. By
    Corollary \ref{cor:ramified-prime-obstruction-suffices}, all other torsion obstruction sets
    associated to $\Gamma$ have nontrivial intersection with
    $\Delta_{\fa,S}$.
\end{ex}

The definite class number one quaternion algebras $B/\bQ$ were
classified by Brzezinski \cite{Brzezinski95}; their
discriminants are
\[
    D(B)=2,3,5,7,13.
\]
We now apply the torsion obstruction criterion to prove that every
uninverted minimal complementary triple associated to one of these
algebras gives rise to infinitely many torsion-free $S$-arithmetic
groups. To ease the notation, we denote all ideals by their generators
in $\bZ$.

\begin{thm}
    Let $B/\bQ$ be a definite quaternion algebra of class number one,
    and let $(T,C,H)$ be an uninverted or trivial minimal triple. Then,
    for every $n \geq 1$, there exist infinitely many $C$-eligible
    finite subsets $S \subseteq V_f(\bQ)$ with $\#S_0 = n$ such that
    $\Gamma_{\fO,S}(C,H)$ is torsion-free.
\end{thm}

\begin{proof}
    Fix a triple $(T,C,H)$. Suppose first that, for some integer $a$, the relevant
    obstruction set does not contain $1$. If every $p \in S_0$ satisfies
    \[
        p \equiv 1 \bmod a,
    \]
    then $\Delta_{a,S} = 1$, and Proposition
    \ref{prop:torsion-obstruction-criterion} implies that
    $\Gamma_{\fO,S}(C,H)$ is torsion-free. By Dirichlet's theorem on
    arithmetic progressions \cite[Exercise I.10.1, Theorem VII.5.14]{Neukirch99}, there are infinitely many primes
    $p \equiv 1 \bmod a$. Hence, for any prescribed cardinality $n$, there
    are infinitely many choices of $S$ with $\#S_0=n$.

    For every nontrivial triple $(T,C,H) \in \cM(B)$, Appendix
    \ref{app:complementary-triple-tables} provides one of
    \[
        \Theta_{\fe(C)}(\widetilde H)
        \qquad \text{or} \qquad
        \Theta_{\fe(D)}
    \]
    that does not contain $1$. The only remaining cases are the trivial
    triples for the maximal order $\fO \subseteq B$ with $D(B)=13$.
    Appendix \ref{app:complementary-triple-tables} gives the same
    conclusion in these cases, since $1 \notin \Theta_{13}$.
\end{proof}

\begin{rmk}
    A closely related trace-zero criterion appears in work of
    Amorós--Milione \cite{AmorosMilione18} on $p$-adic uniformization of
    Shimura curves. They consider a class number one Eichler order
    $\fO$ in a definite rational quaternion algebra and principal
    congruence subgroups
    \[
        \Gamma_p(\xi)=
        \left\{\alpha\in \fO[1/p]^\times:\alpha\equiv 1 \bmod \xi \fO\right\}/
        \bZ[1/p]^\times, \quad \xi \in \fO.
    \]
    Under their right-unit hypothesis on $\xi$, the group
    $\Gamma_p(\xi)$ is generated by the classes of elements
    $\alpha\in \fO$ with $\Nrd(\alpha)=p$ and
    $\alpha\equiv 1\bmod \xi \fO$. The trace-zero generators are
    precisely the involutions in this generating set, and their
    null-trace condition
    \[
        \#\left\{\alpha\in \fO:\Nrd(\alpha)=p,\ \alpha\equiv 1\bmod \xi \fO,\ 
        \trd(\alpha)=0\right\}=0
    \]
    implies that $\Gamma_p(\xi)$ is a Schottky group. It would be
    interesting to extend their results to other totally definite
    orders, but we do not pursue this here.
\end{rmk}

\subsection{The Hurwitz order} \label{sec:hurwitz-order}

We now give an extended example with the well-studied algebra
$B = \quat{\bQ}{-1}{-1}$. A maximal order in $B$ is the
\emph{Hurwitz order}
\[
    \fO = \bZ\cdot 1 + \bZ\cdot i + \bZ\cdot k
        + \bZ\cdot \frac{1+i+j+k}{2}.
\]
Reduction modulo $2$ gives an isomorphism
\[
    \pi_{\fO,2}:\Gamma_{\fO,\emptyset} \to G_\fO[2].
\]
Hence, we have a complementary triple $(\emptyset,2,\{1\})$.
Furthermore, we find
\[
    \Theta_4(\widetilde H) = \left\{-1\right\}.
\]
Therefore, if $S$ is a finite set of odd primes $p \equiv 1 \bmod 4$,
then $\Gamma_{\fO,S}(2)$ is torsion-free. However, we can find more
torsion-free families by considering complementary triples with $C = 4$.

Indeed, we find three distinct conjugacy classes of complements to the
image of $\Gamma_{\fO,\emptyset}$ in $G_\fO[4]$. The corresponding
groups are
\[
    H_1 \cong \bZ/2\bZ \times \bZ/4\bZ,
    \quad
    H_2 \cong (\bZ/2\bZ)^3,
    \quad
    H_3 \cong D_4.
\]
Moreover,
\[
    \Theta_4(H_1) = \emptyset,
    \quad
    \Theta_4(H_2) = \left\{3\right\},
    \quad
    \Theta_4(H_3) = \left\{1,3\right\}.
\]
Therefore, the groups $\Gamma_{\fO,S}(4,H_1)$ are always torsion-free
by Corollary \ref{cor:empty-obstruction-set-torsion-free}. On the other hand, we have
\[
    H_2 = \ker(\pi_{4,2}),
\]
where $\pi_{4,2}:G_\fO[4]\to G_\fO[2]$ is reduction modulo $2$. Hence,
for all eligible $S$, there is an equality
\[
    \Gamma_{\fO,S}(4,H_2) = \Gamma_{\fO,S}(2).
\]

We now give an alternative construction of the groups
$\Gamma_{\fO,S}(4,H_1)$ and $\Gamma_{\fO,S}(4,H_2)$. Inside $\fO$ is
the \emph{Lipschitz order}
\[
    \fL = \bZ\cdot 1 + \bZ\cdot i + \bZ\cdot j + \bZ \cdot k,
\]
which is not Eichler. However,
\[
    \fO[1/2] = \fL[1/2],
\]
and $\fL_p$ is maximal for every odd prime $p$. Since elements of
$\fL$ are especially easy to write down, we do so in this case. We have
\[
    \Gamma_{\fL,\emptyset}
    =
    \left\{1,i,j,ij\right\}
    \cong
    (\bZ/2\bZ)^2,
\]
and
\[
    G_\fL[2]
    =
    \left\{
    \begin{array}{cccc}
    1, & j, & 1+i+j, & 1+i+k, \\
    i, & k, & 1+j+k, & i+j+k
    \end{array}
    \right\}
    \cong
    (\bZ/2\bZ)^3.
\]
It follows immediately that the map
\[
    \Gamma_{\fL,\emptyset} \to G_\fL[2]
\]
is injective. Furthermore, there are four conjugacy classes of
complements, given by
\[
    K_1 = \langle i+j+k \rangle,
    \quad
    K_2 = \langle 1+i+j \rangle,
    \quad
    K_3 = \langle 1+i+ij \rangle,
    \quad
    K_4 = \langle 1+j+ij \rangle.
\]
In particular, the triples $(\emptyset,2,K_l)$ yield torsion-free
groups for $l=2,3,4$, since these subgroups $K_l$ contain no elements
of trace zero.

We now show how the groups $\Gamma_{\fL,S}(2,K_l)$ embed into
$\Gamma_{\fO,S}$. First, observe from the given bases that
$\Gamma_{\fO,S}(2) \subseteq \Gamma_{\fL,S}$. For every $k \geq 1$, the
restriction of $\pi_{2^k}$ to $\Gamma_{\fL,S}$ fits into the following
commutative diagram:
\[\begin{tikzcd}
	{\Gamma_{\fL,S}} & {\Gamma_{\fO,S}} \\
	{G_\fL[2^k]} & {G_\fO[2^k]}
	\arrow[tail, from=1-1, to=1-2]
	\arrow["{\pi_{2^k}}"', from=1-1, to=2-1]
	\arrow["{\pi_{2^k}}", from=1-2, to=2-2]
	\arrow["{\iota_{k}}"', from=2-1, to=2-2]
\end{tikzcd}\]
The bottom map $\iota_k$ need not be injective or surjective. For
example, $\iota_1$ has kernel $K_1$, and
\[
    \Gamma_{\fO,S}(2) = \Gamma_{\fL,S}(2,K_1) \leq \Gamma_{\fL,S}
\]
for finite sets of odd primes $S$. For $l=2,3,4$, there is a unique
lift, up to conjugacy, of $K_l$ to $G_\fL[4]$. These lifted groups all
map under $\iota_2$ to $H_1$. Therefore,
\[
    \Gamma_{\fO,S}(4,H_1) = \Gamma_{\fL,S}(2,K_l) \leq \Gamma_{\fL,S}
\]
for $l = 2,3,4$.

The above discussion about the Hurwitz and Lipschitz orders recovers
results of Lubotzky--Phillips--Sarnak \cite{LPS88},
Rattaggi \cite{Rattaggi04}, and Rungtanapirom--Stix--Vdovina
\cite{RungtanapiromStixVdovina19}.

\subsection{Irreducible complementary triples} \label{sec:irreducible-complementary-triples}

We call a triple $(T,\fC,H)$ \emph{irreducible} if, for all
$\fA \mid \fC$ with $\fA \neq \fC$ such that $\pi_\fA$ is injective on
$\Gamma_T$, the group $\pi_{\fC,\fA}(H)$ is not a complement of
$\pi_\fA(\Gamma_T)$ in $G[\fA]$. Clearly, a minimal triple is
irreducible, as $\pi_\fA$ is not injective on $\Gamma_T$ for all
$\fA \mid \fC$ with $\fA \neq \fC$. However, as we saw in the previous
section, the converse is not true, and irreducible, non-minimal triples
can have a rich structure.

The most interesting irreducible, non-minimal complementary triples
$(T,\fC,H)$ are those with $\fC = \fp^e$ for a prime $\fp$ and
\[
    e > e_0(\fp,T),
    \qquad
    e_0(\fp,T)
    =
    \min\left\{e \geq 1 : \pi_{\fp^e}
    \text{ is injective on $\Gamma_T$}\right\},
\]
but
\[
    H_0 \cap \pi_{\fp^{e_0}}(\Gamma_T) \neq 1,
    \qquad
    H_0 = \pi_{\fp^e,\fp^{e_0}}(H),
\]
where $e_0=e_0(\fp,T)$. We saw such examples in the previous section,
and we can give restrictions on when such triples occur.

\begin{prop}
    Let $\fp$ be a prime ideal with residue characteristic $p$, and
    suppose that $p \nmid \#\Gamma_T$. Let $e \geq e_0=e_0(\fp,T)$,
    and suppose that $(T,\fp^e,H)$ is a complementary triple of $\fO$.
    If
    \[
        H_0 = \pi_{\fp^e,\fp^{e_0}}(H),
    \]
    then $(T,\fp^{e_0},H_0)$ is a complementary triple of $\fO$. In
    particular, if $e>e_0$, then $(T,\fp^e,H)$ is reducible.
\end{prop}

\begin{proof}
    Let $K_0 = \ker(\pi_{\fp^e,\fp^{e_0}})$. Given
    $\gamma \in \Gamma_T \cap H_0$, it follows from the triviality of
    $\Gamma_T \cap H$ that there exists $h_\gamma \in H$ such that
    $\gamma h_\gamma^{-1} \in K_0$. We get a bijection
    \[
        H_0 \cap \Gamma_T
        \to
        K_0/(H \cap K_0),
        \quad
        \gamma \mapsto \gamma h_\gamma^{-1}(H \cap K_0),
    \]
    so $\#(H_0 \cap \Gamma_T) \mid \#K_0$. In both the ramified and split
    cases, $K_0$ is a $p$-group. Thus, $H_0 \cap \Gamma_T$ is a $p$-group
    if it is nontrivial. Since $p\nmid \#\Gamma_T$, we must have
    $H_0 \cap \Gamma_T=1$, so $(T,\fp^{e_0},H_0)$ is a complementary
    triple.
\end{proof}

The proposition shows that prime-power levels rarely produce new
irreducible examples. At composite levels, however, non-minimal
irreducible triples can occur in infinite families. The obstacle is
that a complementary subgroup modulo $\fA\fB$ need not project to a
complementary subgroup modulo either factor. More precisely, if
$\fA,\fB$ are coprime, nontrivial ideals and $\pi_\fA$ is injective on
$\Gamma_T$, there can exist complementary triples $(T,\fA\fB,H)$ such
that neither $\pi_\fA(H)$ nor $\pi_\fB(H)$ is complementary to the
corresponding image of $\Gamma_T$.

For example, let $\fO$ be the Hurwitz order and let $T=\{2\}$. Then
$(1+i)\cdot\fO$ is the unique two-sided ideal above $2$, and
\[
    \lag [1+i] \rag \leq \wtU_2^1
\]
is isomorphic to $C_2$. Let $p$ be an odd prime such that $2$ is not a
square modulo $p$. Under the isomorphism
\[
    G[2p]\cong G[2]\times G[p],
\]
set
\[
    H=\lag [1+i] \rag \times \PSL_2(\bF_p).
\]
We show that $H$ is complementary to $\Gamma_T$. First, observe that
\[
    \#H\cdot \#\Gamma_T
    =
    2\cdot \frac{\#G[p]}{2}\cdot \#\Gamma_T
    =
    \#G[p]\cdot \#G[2]
    =
    \#G[2p],
\]
where we use $\#\Gamma_T=\#G[2]$. Moreover,
$H\cap \Gamma_T=1$. Indeed, any nontrivial element in the intersection
would have $G[2]$-component $\pi_2(1+i)$, and hence would have
$G[p]$-component $\pi_p(1+i)$. But
$\pi_p(1+i)\notin \PSL_2(\bF_p)$, since $\Nrd(1+i)=2$ is not a square
modulo $p$. Thus, we obtain infinitely many non-minimal, irreducible
triples $(T,2p,H)$.

There is still a useful elementary constraint. If $(T,\fC,H)$ is a
complementary triple and $\fA\mid \fC$, then projecting modulo $\fA$
gives
\[
    \pi_\fA(H)\cdot \pi_\fA(\Gamma_T)=G[\fA],
\]
and hence
\[
    [G[\fA]:\pi_\fA(H)] \leq \#\pi_\fA(\Gamma_T).
\]
In particular, if $\fA=\fp^e$ is a prime power and
\[
    [G[\fA]:\pi_\fA(H)] > 2,
\]
then Dickson's theorem bounds $N\fp$. This suggests that the
classification of irreducible complementary triples may be extendable
beyond the minimal case. We do not pursue this here.

\section{Acknowledgments}

The authors thank the Einstein Institute of Mathematics at the Hebrew University of Jerusalem for hosting them during the summer 2022 REU, where this project began. 
They are grateful to Ari Shnidman for organizing the REU, 
and to their REU advisor, Shai Evra, for introducing them to Ramanujan graphs, quaternion algebras, 
and the central problem addressed in this paper. 
They also thank Professor Evra for reading early drafts of this manuscript and providing inspiration and helpful feedback.

The first author thanks Alina Vdovina for introducing him to her work on simply transitive actions on products of trees, 
and for many subsequent discussions on this topic.
This paper owes a substantial intellectual debt to the work of Vdovina and her co-authors. 

The first author thanks his advisor, Alan Reid, for many useful discussions about arithmetic groups and quaternion algebras, and for his interest in and encouragement of this project. 
The first author also thanks Jakob Stix for a helpful discussion on $S$-arithmetic groups.

The authors also thank Alina Vdovina and Alan Reid for providing feedback on this manuscript.

\section{AI Disclosure}

The authors used AI-assisted tools to help develop and revise the Magma code
used in this work, edit the language and clarity of exposition, format the
manuscript in \LaTeX{}, and prepare the tables and figures included in the
paper. The authors take full responsibility for the content of the manuscript.

\bibliography{references}

\appendix

\section{Complementary triples table}
\label{app:complementary-triple-tables}
\input{complementary_triples_tables}

\end{document}

%% file: commands.tex
\usepackage[utf8]{inputenc}

\usepackage{amsmath, amssymb, amsthm, amsfonts, amsopn, mathrsfs}

\let\emptyset\varnothing

\usepackage{tikz}
\usepackage{tikz-cd}
\usepackage{quiver}
\usepackage{xcolor}
\usepackage{geometry}
\usepackage[bottom]{footmisc}
\usepackage{array}
\usepackage{longtable}
\usepackage{booktabs}
\usepackage{enumitem}

\usepackage[pdfpagelabels]{hyperref}
\hypersetup{
  colorlinks=true,
  breaklinks=true,
  bookmarksopen=true,
  bookmarksnumbered=true,
  pdfpagemode=UseOutlines,
  plainpages=false
}

\newcommand{\bF}{\mathbb{F}}

\newcommand{\bQ}{\mathbb{Q}}
\newcommand{\bR}{\mathbb{R}}

\newcommand{\bZ}{\mathbb{Z}}

\newcommand{\cG}{\mathscr{G}}

\newcommand{\cM}{\mathscr{M}}

\newcommand{\cR}{\mathscr{R}}

\newcommand{\cT}{\mathscr{T}}

\newcommand{\cX}{\mathscr{X}}

\newcommand{\fA}{\mathfrak{A}}
\newcommand{\fB}{\mathfrak{B}}
\newcommand{\fC}{\mathfrak{C}}
\newcommand{\fD}{\mathfrak{D}}
\newcommand{\fE}{\mathfrak{E}}

\newcommand{\fI}{\mathfrak{I}}
\newcommand{\fJ}{\mathfrak{J}}

\newcommand{\fL}{\mathfrak{L}}

\newcommand{\fN}{\mathfrak{N}}
\newcommand{\fO}{\mathfrak{O}}
\newcommand{\fP}{\mathfrak{P}}
\newcommand{\fQ}{\mathfrak{Q}}

\newcommand{\fa}{\mathfrak{a}}
\newcommand{\fb}{\mathfrak{b}}

\newcommand{\fe}{\mathfrak{e}}

\newcommand{\fo}{\mathfrak{o}}
\newcommand{\fp}{\mathfrak{p}}
\newcommand{\fq}{\mathfrak{q}}

\newcommand{\lag}{\left\langle}
\newcommand{\rag}{\right\rangle}
\newcommand{\quat}[3]{\left(\frac{#2,#3}{#1}\right)}
\newcommand{\wtU}{\widetilde{U}}

\DeclareMathOperator{\Adj}{Adj}
\DeclareMathOperator{\Aut}{Aut}
\DeclareMathOperator{\Cay}{Cay}
\DeclareMathOperator{\characteristic}{char}
\DeclareMathOperator{\Cl}{Cl}
\DeclareMathOperator{\Frob}{Frob}
\DeclareMathOperator{\Gal}{Gal}
\DeclareMathOperator{\GL}{GL}
\DeclareMathOperator{\Nrd}{Nrd}
\DeclareMathOperator{\PGL}{PGL}
\DeclareMathOperator{\Pic}{Pic}
\DeclareMathOperator{\PSL}{PSL}
\DeclareMathOperator{\Ram}{Ram}
\DeclareMathOperator{\SO}{SO}
\DeclareMathOperator{\Sym}{Sym}
\DeclareMathOperator{\tr}{tr}
\DeclareMathOperator{\trd}{trd}

\newtheorem{thm}{Theorem}
\newtheorem{prop}[thm]{Proposition}
\newtheorem{lem}[thm]{Lemma}
\newtheorem{cor}[thm]{Corollary}

\newtheorem{thmABC}{Theorem}

\theoremstyle{definition}

\theoremstyle{remark}
\newtheorem{rmk}[thm]{Remark}
\newtheorem{nota}[thm]{Notation}
\newtheorem{ex}[thm]{Example}

\numberwithin{equation}{section}

%% file: complementary_triples_tables.tex
This appendix records $\cM(B)$ for each totally definite quaternion algebra
$B/K$ of class number one in the Kirschmer--Voight classification. 

For a number field $K$, let $n=[K:\bQ]$ and let $d_K$ be its discriminant.
It is a consequence of the Kirschmer--Voight classification that, for each
algebra $B/K$ appearing in the tables, the pair $(n,d_K)$ determines $K$.

Moreover, if two class number one,
totally definite quaternion algebras $B_1$ and $B_2$ have $D(B_1) = D(B_2)$, then there is an
automorphism $\sigma \in \Aut_{\bQ}(K)$ that carries the ramification set of
$B_1$ to that of $B_2$.  Equivalently, there is a $\bQ$-algebra isomorphism
$B_1 \to B_2$ whose restriction to $K$ is $\sigma$.  The image of a maximal
order under this isomorphism is again maximal, and the isomorphism carries
complementary triples to complementary triples. We list our table conventions below.

\begin{itemize}[leftmargin=*, itemsep=0.65em]

\item
Totally definite, class number one algebras are identified up to the invariants $n$, $d_K$, and $D$.

\item
We record whether $1 \in \Theta_{\fe(\fD)}$ (see Section \ref{sec:torsion-free-groups} for notation).
This condition is invariant under a $\bQ$-isomorphism of
quaternion algebras. Thus, in the cases where the table identifies two algebras, the condition agrees for both. 
When $\fD(B)=\fo_K$, this column contains $-$.

\item
The set $T \subseteq \Ram_f(B)$ consists of prime ideals of $\fo_K$, and the complementary ideal $\fC$ is written
as a product of prime ideals. We use the following convention to write prime
ideals in $\fo_K$:
\begin{enumerate}
    \item $\fp_{q,\mathrm{ram}}$ denotes a prime ideal of norm $q$ ramified
    in $B$,
    \item $\fp_{q,\mathrm{spl}}$ denotes a prime ideal of norm $q$ split
    in $B$.
\end{enumerate}
When distinct prime ideals of the same norm and the same ramification status
occur among the rows with fixed $d_K$, $D$, and $T$, subscripts
$\mathrm{a},\mathrm{b},\mathrm{c},\ldots$ are added.

\item
When $\Gamma_\emptyset$ is nontrivial, no row is included for a subset
$T \subseteq \Ram_f(B)$ unless a nontrivial complementary triple
$(T,\fC,H)$ exists for that $T$.

\item
When $\Gamma_\emptyset=1$, recall from Section~\ref{sec:computing-complementary-triples} that
certain complementary triples $(T,\fC,H)$ with $\#T=0$ or $\#T=1$ are treated
as trivial.  We do not list the triples covered by that convention.  Instead,
the row with $T=\emptyset$ and $-$ in the complementary-triple columns records the existence 
of the algebra. Note that if $\#\Ram_f(B) > 1$, then there may still be triples in $\cM(B)$, and those triples are
listed when they occur.

\item
The columns $c$ and $c^\triangleleft$ denote, respectively, the number of
conjugacy classes of complements $H$ and the number of such conjugacy classes
represented by normal complements.

\item
For every triple $(T,\fC,H)$ such that no prime in $T$ divides $\fC$, we set
$\widetilde H=\pi_{\fe(\fC),\fC}^{-1}(H)$ and compute
\[
    \Theta = \Theta_{\fe(\fC)}(\widetilde H).
\]
The final two columns of the table record $\#\Theta$ and whether this obstruction
set contains $1$. The symbol $\star$ indicates that $1 \in \Theta$; otherwise,
the entry is left blank. A dash in both columns means that no torsion obstruction set is recorded for
that row; this occurs, in particular, when some prime in $T$ divides $\fC$.

\item
In the rare cases where distinct conjugacy classes of complements for the same
$T$ and $\fC$ yield different values of $(\#\Theta,\,1 \in \Theta)$, the
corresponding pairs are listed on separate lines.

\end{itemize}

\begingroup
\setlength{\tabcolsep}{3.5pt}
\renewcommand{\arraystretch}{1.18}

\begin{longtable}{@{}r r r@{\hspace{0.45em}}c@{\hspace{0.65em}}l l@{\hspace{0.55em}}>{\raggedright\arraybackslash}p{0.19\textwidth}@{\hspace{0.6em}}r@{\hspace{0.6em}}r@{\hspace{0.6em}}c@{\hspace{0.6em}}c@{}}
\caption{Complementary triples.}\label{tab:complementary-triples}\\
\toprule
$n$ & $d_K$ & $D$ & $1\in\Theta_{\fe(\fD)}$ & $T$ & $\Gamma_T$ & $\fC$ & $c$ & $c^\triangleleft$ & $\#\Theta$ & $1\in\Theta$ \\
\midrule
\endfirsthead
\toprule
$n$ & $d_K$ & $D$ & $1\in\Theta_{\fe(\fD)}$ & $T$ & $\Gamma_T$ & $\fC$ & $c$ & $c^\triangleleft$ & $\#\Theta$ & $1\in\Theta$ \\
\midrule
\endhead
1 & 1 & 2 & $\star$ & $\emptyset$ & $A_4$ & $\fp_{2,\mathrm{ram}}$ & 1 & 1 & 1 &  \\
 &  &  &  &  &  & $\fp_{3,\mathrm{spl}}$ & 1 & 0 & 1 &  \\
 &  &  &  &  &  & $\fp_{11,\mathrm{spl}}$ & 1 & 0 & 1 &  \\
\addlinespace[3pt]
 &  &  &  & $\{\fp_{2,\mathrm{ram}}\}$ & $S_4$ & $\fp_{2,\mathrm{ram}}$ & 1 & 1 & $-$ & $-$ \\
 &  &  &  &  &  & $\fp_{3,\mathrm{spl}}$ & 1 & 1 & 0 &  \\
 &  &  &  &  &  & $\fp_{5,\mathrm{spl}}$ & 1 & 0 & 0 &  \\
 &  &  &  &  &  & $\fp_{7,\mathrm{spl}}$ & 1 & 0 & 1 &  \\
 &  &  &  &  &  & $\fp_{11,\mathrm{spl}}$ & 1 & 0 & 0 &  \\
 &  &  &  &  &  & $\fp_{23,\mathrm{spl}}$ & 1 & 0 & 1 &  \\
\midrule
 &  & 3 & $\star$ & $\emptyset$ & $S_3$ & $\fp_{2,\mathrm{spl}}$ & 1 & 1 & 2 &  \\
\addlinespace[3pt]
 &  &  &  & $\{\fp_{3,\mathrm{ram}}\}$ & $D_6$ & $\fp_{3,\mathrm{ram}}$ & 2 & 0 & $-$ & $-$ \\
 &  &  &  &  &  & $\fp_{2,\mathrm{spl}}^{2}$ & 4 & 1 & \begin{tabular}[t]{@{}c@{}}0\\ 1\end{tabular} &  \\
 &  &  &  &  &  & $\fp_{11,\mathrm{spl}}$ & 1 & 0 & 1 &  \\
\midrule
 &  & 5 &  & $\emptyset$ & $C_3$ & $\fp_{2,\mathrm{spl}}$ & 1 & 0 & 4 & $\star$ \\
 &  &  &  &  &  & $\fp_{3,\mathrm{spl}}$ & 1 & 0 & 2 & $\star$ \\
 &  &  &  &  &  & $\fp_{5,\mathrm{ram}}$ & 1 & 1 & 1 &  \\
\addlinespace[3pt]
 &  &  &  & $\{\fp_{5,\mathrm{ram}}\}$ & $S_3$ & $\fp_{2,\mathrm{spl}}$ & 1 & 1 & 2 &  \\
 &  &  &  &  &  & $\fp_{3,\mathrm{spl}}$ & 2 & 1 & 1 & $\star$ \\
 &  &  &  &  &  & $\fp_{5,\mathrm{ram}}$ & 2 & 1 & $-$ & $-$ \\
\midrule
 &  & 7 & $\star$ & $\emptyset$ & $C_2$ & $\fp_{2,\mathrm{spl}}$ & 1 & 1 & 2 &  \\
\addlinespace[3pt]
 &  &  &  & $\{\fp_{7,\mathrm{ram}}\}$ & $C_2^2$ & $\fp_{3,\mathrm{spl}}$ & 1 & 0 & 1 &  \\
 &  &  &  &  &  & $\fp_{2,\mathrm{spl}}^{2}$ & 1 & 1 & 0 &  \\
\midrule
 &  & 13 &  & $\emptyset$ & $1$ & $-$ & $-$ & $-$ & $-$ & $-$ \\
\midrule
\midrule
2 & 5 & 1 & $-$ & $\emptyset$ & $A_5$ & $\fp_{4,\mathrm{spl}}$ & 1 & 1 & 2 &  \\
 &  &  &  &  &  & $\fp_{5,\mathrm{spl}}$ & 1 & 0 & 1 &  \\
 &  &  &  &  &  & $\fp_{11,\mathrm{spl},\mathrm{a}}$ & 1 & 0 & 1 &  \\
 &  &  &  &  &  & $\fp_{11,\mathrm{spl},\mathrm{b}}$ & 1 & 0 & 1 &  \\
 &  &  &  &  &  & $\fp_{59,\mathrm{spl},\mathrm{a}}$ & 1 & 0 & 1 &  \\
 &  &  &  &  &  & $\fp_{59,\mathrm{spl},\mathrm{b}}$ & 1 & 0 & 1 &  \\
\midrule
 &  & 20 &  & $\emptyset$ & $C_5$ & $\fp_{4,\mathrm{ram}}$ & 1 & 1 & 7 & $\star$ \\
\addlinespace[3pt]
 &  &  &  & $\{\fp_{4,\mathrm{ram}}\}$ & $D_5$ & $\fp_{4,\mathrm{ram}}$ & 4 & 1 & $-$ & $-$ \\
\addlinespace[3pt]
 &  &  &  & $\{\fp_{5,\mathrm{ram}}\}$ & $C_{10}$ & $\fp_{9,\mathrm{spl}}$ & 1 & 0 & 1 & $\star$ \\
\midrule
 &  & 44 & $\star$ & $\{\fp_{4,\mathrm{ram}}\}$ & $C_4$ & $\fp_{11,\mathrm{ram}}$ & 1 & 1 & 0 &  \\
\addlinespace[3pt]
 &  &  &  & $\{\fp_{11,\mathrm{ram}}\}$ & $C_2^2$ & $\fp_{11,\mathrm{ram}}$ & 1 & 0 & $-$ & $-$ \\
\addlinespace[3pt]
 &  &  &  & $\{\fp_{4,\mathrm{ram}},\allowbreak \fp_{11,\mathrm{ram}}\}$ & $D_4$ & $\fp_{11,\mathrm{ram}}$ & 1 & 1 & $-$ & $-$ \\
\midrule
\midrule
2 & 8 & 1 & $-$ & $\emptyset$ & $S_4$ & $\fp_{2,\mathrm{spl}}^{2}$ & 3 & 1 & \begin{tabular}[t]{@{}c@{}}2\\ 6\end{tabular} &  \\
 &  &  &  &  &  & $\fp_{7,\mathrm{spl},\mathrm{a}}$ & 1 & 0 & 1 &  \\
 &  &  &  &  &  & $\fp_{7,\mathrm{spl},\mathrm{b}}$ & 1 & 0 & 1 &  \\
 &  &  &  &  &  & $\fp_{23,\mathrm{spl},\mathrm{a}}$ & 1 & 0 & 1 &  \\
 &  &  &  &  &  & $\fp_{23,\mathrm{spl},\mathrm{b}}$ & 1 & 0 & 1 &  \\
\midrule
 &  & 14 & $\star$ & $\{\fp_{2,\mathrm{ram}}\}$ & $C_8$ & $\fp_{7,\mathrm{ram}}$ & 1 & 1 & 0 &  \\
 &  &  &  &  &  & $\fp_{7,\mathrm{spl}}$ & 1 & 0 & 1 &  \\
\addlinespace[3pt]
 &  &  &  & $\{\fp_{7,\mathrm{ram}}\}$ & $D_4$ & $\fp_{7,\mathrm{ram}}$ & 1 & 0 & $-$ & $-$ \\
\addlinespace[3pt]
 &  &  &  & $\{\fp_{2,\mathrm{ram}},\allowbreak \fp_{7,\mathrm{ram}}\}$ & $D_8$ & $\fp_{7,\mathrm{ram}}$ & 1 & 1 & $-$ & $-$ \\
 &  &  &  &  &  & $\fp_{7,\mathrm{spl}}$ & 1 & 0 & 0 &  \\
\midrule
 &  & 18 &  & $\emptyset$ & $C_3$ & $\fp_{2,\mathrm{ram}}$ & 1 & 1 & 7 & $\star$ \\
\addlinespace[3pt]
 &  &  &  & $\{\fp_{2,\mathrm{ram}}\}$ & $S_3$ & $\fp_{2,\mathrm{ram}}$ & 2 & 1 & $-$ & $-$ \\
\midrule
 &  & 50 &  & $\emptyset$ & $1$ & $-$ & $-$ & $-$ & $-$ & $-$ \\
\addlinespace[3pt]
 &  &  &  & $\{\fp_{2,\mathrm{ram}}\}$ & $C_2$ & $\fp_{25,\mathrm{ram}}$ & 1 & 1 & 0 &  \\
\addlinespace[3pt]
 &  &  &  & $\{\fp_{2,\mathrm{ram}},\allowbreak \fp_{25,\mathrm{ram}}\}$ & $C_2^2$ & $\fp_{25,\mathrm{ram}}$ & 1 & 1 & $-$ & $-$ \\
\midrule
\midrule
2 & 13 & 1 & $-$ & $\emptyset$ & $A_4$ & $\fp_{3,\mathrm{spl},\mathrm{a}}$ & 1 & 0 & 1 &  \\
 &  &  &  &  &  & $\fp_{3,\mathrm{spl},\mathrm{b}}$ & 1 & 0 & 1 &  \\
 &  &  &  &  &  & $\fp_{4,\mathrm{spl}}$ & 1 & 0 & 2 &  \\
\midrule
 &  & 12 & $\star$ & $\{\fp_{3,\mathrm{ram}}\}$ & $C_2^2$ & $\fp_{3,\mathrm{ram}}$ & 1 & 0 & $-$ & $-$ \\
\addlinespace[3pt]
 &  &  &  & $\{\fp_{4,\mathrm{ram}}\}$ & $C_4$ & $\fp_{3,\mathrm{ram}}$ & 1 & 1 & 0 &  \\
 &  &  &  &  &  & $\fp_{3,\mathrm{spl}}$ & 1 & 0 & 1 &  \\
\addlinespace[3pt]
 &  &  &  & $\{\fp_{3,\mathrm{ram}},\allowbreak \fp_{4,\mathrm{ram}}\}$ & $D_4$ & $\fp_{3,\mathrm{ram}}$ & 1 & 1 & $-$ & $-$ \\
 &  &  &  &  &  & $\fp_{3,\mathrm{spl}}$ & 1 & 0 & 0 &  \\
\midrule
\midrule
2 & 17 & 1 & $-$ & $\emptyset$ & $S_3$ & $\fp_{2,\mathrm{spl},\mathrm{a}}$ & 1 & 1 & 2 &  \\
 &  &  &  &  &  & $\fp_{2,\mathrm{spl},\mathrm{b}}$ & 1 & 1 & 2 &  \\
\midrule
\midrule
3 & 49 & 7 & $\star$ & $\{\fp_{7,\mathrm{ram}}\}$ & $D_{14}$ & $\fp_{27,\mathrm{spl}}$ & 1 & 0 & 1 &  \\
\midrule
 &  & 8 & $\star$ & $\{\fp_{8,\mathrm{ram}}\}$ & $S_4$ & $\fp_{7,\mathrm{spl}}$ & 1 & 0 & 1 &  \\
\midrule
 &  & 13 &  & $\emptyset$ & $C_7$ & $\fp_{13,\mathrm{ram}}$ & 1 & 1 & 1 &  \\
\addlinespace[3pt]
 &  &  &  & $\{\fp_{13,\mathrm{ram}}\}$ & $D_7$ & $\fp_{7,\mathrm{spl}}$ & 1 & 0 & 1 & $\star$ \\
 &  &  &  &  &  & $\fp_{13,\mathrm{ram}}$ & 2 & 1 & $-$ & $-$ \\
 &  &  &  &  &  & $\fp_{13,\mathrm{spl}}$ & 1 & 0 & 1 & $\star$ \\
\midrule
 &  & 29 &  & $\emptyset$ & $C_3$ & $\fp_{29,\mathrm{ram}}$ & 1 & 1 & 1 &  \\
\addlinespace[3pt]
 &  &  &  & $\{\fp_{29,\mathrm{ram}}\}$ & $S_3$ & $\fp_{29,\mathrm{ram}}$ & 2 & 1 & $-$ & $-$ \\
\midrule
 &  & 43 & $\star$ & $\{\fp_{43,\mathrm{ram}}\}$ & $C_2^2$ & $\fp_{43,\mathrm{ram}}$ & 1 & 0 & $-$ & $-$ \\
\midrule
\midrule
3 & 81 & 3 & $\star$ & $\{\fp_{3,\mathrm{ram}}\}$ & $D_{18}$ & $\fp_{3,\mathrm{ram}}\fp_{8,\mathrm{spl}}$ & 1 & 0 & $-$ & $-$ \\
\midrule
 &  & 19 & $\star$ & $\{\fp_{19,\mathrm{ram}}\}$ & $C_2^2$ & $\fp_{3,\mathrm{spl}}$ & 1 & 0 & 1 &  \\
 &  &  &  &  &  & $\fp_{19,\mathrm{ram}}$ & 1 & 0 & $-$ & $-$ \\
\midrule
 &  & 37 &  & $\emptyset$ & $1$ & $-$ & $-$ & $-$ & $-$ & $-$ \\
\midrule
\midrule
3 & 148 & 2 & $\star$ & $\{\fp_{2,\mathrm{ram}}\}$ & $S_4$ & $\fp_{5,\mathrm{spl}}$ & 1 & 0 & 0 &  \\
 &  &  &  &  &  & $\fp_{23,\mathrm{spl}}$ & 1 & 0 & 1 &  \\
\midrule
 &  & 5 &  & $\emptyset$ & $C_3$ & $\fp_{2,\mathrm{spl}}$ & 1 & 0 & 16 & $\star$ \\
 &  &  &  &  &  & $\fp_{5,\mathrm{ram}}$ & 1 & 1 & 1 &  \\
\addlinespace[3pt]
 &  &  &  & $\{\fp_{5,\mathrm{ram}}\}$ & $S_3$ & $\fp_{2,\mathrm{spl}}$ & 1 & 1 & 8 & $\star$ \\
 &  &  &  &  &  & $\fp_{5,\mathrm{ram}}$ & 2 & 1 & $-$ & $-$ \\
\midrule
 &  & 13 &  & $\emptyset$ & $1$ & $-$ & $-$ & $-$ & $-$ & $-$ \\
\midrule
\midrule
3 & 169 & 5 &  & $\emptyset$ & $C_3$ & $\fp_{5,\mathrm{ram}}$ & 1 & 1 & 1 &  \\
\addlinespace[3pt]
 &  &  &  & $\{\fp_{5,\mathrm{ram}}\}$ & $S_3$ & $\fp_{5,\mathrm{ram}}$ & 2 & 1 & $-$ & $-$ \\
 &  &  &  &  &  & $\fp_{5,\mathrm{spl}}$ & 1 & 0 & 1 & $\star$ \\
\midrule
 &  & 13 &  & $\emptyset$ & $1$ & $-$ & $-$ & $-$ & $-$ & $-$ \\
\midrule
\midrule
3 & 316 & 2 &  & $\emptyset$ & $C_3$ & $\fp_{2,\mathrm{ram}}$ & 1 & 1 & 7 &  \\
 &  &  &  &  &  & $\fp_{2,\mathrm{spl}}$ & 1 & 0 & 4 & $\star$ \\
\addlinespace[3pt]
 &  &  &  & $\{\fp_{2,\mathrm{ram}}\}$ & $S_3$ & $\fp_{2,\mathrm{ram}}$ & 2 & 1 & $-$ & $-$ \\
 &  &  &  &  &  & $\fp_{2,\mathrm{spl}}$ & 1 & 1 & 2 &  \\
\midrule
\midrule
3 & 321 & 3 & $\star$ & $\{\fp_{3,\mathrm{ram}}\}$ & $C_2^2$ & $\fp_{3,\mathrm{ram}}$ & 1 & 0 & $-$ & $-$ \\
 &  &  &  &  &  & $\fp_{3,\mathrm{spl}}$ & 1 & 0 & 1 &  \\
\midrule
\midrule
4 & 725 & 1 & $-$ & $\emptyset$ & $A_5$ & $\fp_{11,\mathrm{spl},\mathrm{a}}$ & 1 & 0 & 1 &  \\
 &  &  &  &  &  & $\fp_{11,\mathrm{spl},\mathrm{b}}$ & 1 & 0 & 1 &  \\
\midrule
\midrule
4 & 1957 & 1 & $-$ & $\emptyset$ & $A_4$ & $\fp_{3,\mathrm{spl}}$ & 1 & 0 & 1 &  \\
\midrule
\midrule
4 & 2777 & 1 & $-$ & $\emptyset$ & $S_3$ & $\fp_{2,\mathrm{spl}}$ & 1 & 1 & 2 &  \\
\midrule
\midrule
5 & 24217 & 5 &  & $\emptyset$ & $C_3$ & $\fp_{5,\mathrm{ram}}$ & 1 & 1 & 1 &  \\
\addlinespace[3pt]
 &  &  &  & $\{\fp_{5,\mathrm{ram}}\}$ & $S_3$ & $\fp_{5,\mathrm{ram}}$ & 2 & 1 & $-$ & $-$ \\
\bottomrule
\end{longtable}

\endgroup